\documentclass{amsart}
\usepackage{amsmath, amsthm}
\usepackage{amsfonts,amssymb}
\usepackage{epsfig}
\usepackage{graphicx}
\usepackage{subfigure}

\newtheorem{theorem}{Theorem}[section]
\newtheorem{cor}{Corollary}
\newtheorem{lemma}[theorem]{Lemma}
\newtheorem{prop}[theorem]{Proposition}

\theoremstyle{definition}
\newtheorem{definition}[theorem]{Definition}
\newtheorem{example}[theorem]{Example}
\newtheorem{remark}[theorem]{Remark}

\graphicspath{{D:/TEXpapers/Continuations/}}
\input{tcilatex}

\begin{document}
\title[Fractal continuation]{Fractal Continuation}
\author{Michael F. Barnsley }
\address{Department of Mathematics\\
Australian National University\\
Canberra, ACT, Australia}
\email{michael.barnsley@maths.anu.edu.au, mbarnsley@aol.com}
\urladdr{http://www.superfractals.com}
\author{Andrew Vince}
\address{Department of Mathematics\\
University of Florida\\
Gainesville, FL 32611-8105, USA}
\email{avince@ufl.edu}

\begin{abstract}
A fractal function is a function whose graph is the attractor of an iterated
function system. This paper generalizes analytic continuation of an analytic
function to continuation of a fractal function.
\end{abstract}

\maketitle

\section{\label{introsec} Introduction}

Analytic continuation is a central concept of mathematics. Riemannian
geometry emerged from the continuation of real analytic functions. This
paper generalizes analytic continuation of an analytic function to
continuation of a fractal function. By fractal function, we mean basically a
function whose graph is the attractor of an iterated function system. We
demonstrate how analytic continuation of a function, defined locally by
means of a Taylor series expansion, generalises to continuation of a, not
necessarily analytic, fractal function.

Fractal functions have a long history, see \cite{tricot} and \cite[Chapter 5]%
{massopust}. They were introduced, in the form considered here, in \cite%
{finterp}. They include many well-known types of non-differentiable
functions, including Takagi curves, Kiesswetter curves, Koch curves,
space-filling curves, and nowhere differentiable functions of Weierstrass. A
fractal function is a continuous function that maps a compact interval $%
I\subset \mathbb{R}$ into a complete metric space, usually $\mathbb{R}$ or $%
\mathbb{R}^{2}$, and may interpolate specified data and have specified
non-integer Minkowski dimension. Fractal functions are the basis of a
constructive approximation theory for non-differentiable functions. They
have been developed both in theory and applications by many authors, see for
example \cite{finterp2, massopust0, massopust, navascues, prasad, scealy,
tosan} and references therein.

Let $N$ be an integer and $\mathcal{I=\{}1,2,...,N\}$. Let $M\geq 2$ be an
integer and $\mathbb{X\subset R}^{M}$ complete with respect to a metric $d_{%
\mathbb{X}}$ that induces, on $\mathbb{X}$, the same topology as the
Euclidean metric. Let $\mathcal{W}$ be an iterated function system (IFS) of
the form 
\begin{equation}  \label{form}
\mathcal{W}=\{\mathbb{X};w_{n},n\in \mathcal{I\}}\text{, }
\end{equation}%
We say that $\mathcal{W}$ is an \textbf{analytic IFS} if $w_{n}$ is a
homeomorphism from ${\mathbb{X}}$ onto ${\mathbb{X}}$ for all $n\in \mathcal{%
I}$, and $w_{n}$ and its inverse $w_{n}^{-1}$ are analytic. By $w_{n}$
analytic, we mean that 
\begin{equation*}
w_{n}(x)=(w_{n1}(x),w_{n2}(x),...,w_{nM}(x)),
\end{equation*}%
where each real-valued function $w_{nm}(x)=w_{nm}(x_{1},x_{2},...,x_{M})$ is
infinitely differentiable in $x_{i}$ with $x_{j}$ fixed for all $j\neq i$,
with a convergent multivariable Taylor series expansion convergent in a
neigbourhood of each point $(x_{1},x_{2},...,x_{M})\in \mathbb{X}$.

To introduce the main ideas, define a \textbf{fractal function} as a
continuous function $f:I\rightarrow\mathbb{R}^{M-1}$, where $I\subset\mathbb{%
R}$ is a compact interval whose graph $G(f)$ is the attractor of a IFS for
the form in eqaution~\eqref{form}. A slightly more restrictive definition
will be given in Section~\ref{fifsec}. If $\mathcal{W}$ is an analytic IFS,
then $f$ is called an \textbf{analytic fractal function}.

The adjective ``fractal" is used to emphasize that $G(f)$ may have
noninteger Hausdorff and Minkowski dimensions. But $f$ may be many times
differentiable or $f$ may even be a real analytic function. Indeed, we prove
that all real analytic functions are, locally, analytic fractal functions;
see Theorem \ref{corollary}. An alternative name for a fractal function $f$
could be a ``self-similar function" because $G(f)$ is a union of transformed
``copies" of itself, specifically 
\begin{equation}  \label{eq:ss}
G(f)=\bigcup_{n=1}^{N}w_{n}(G(f)).
\end{equation}

The goal of this paper is to introduce a new method of analytic
continuation, a method that applies to fractal functions as well as analytic
functions. We call this method \textit{fractal continuation}. When fractal
continuation is applied to a locally defined real analytic function, it
yields the standard analytic continuation. When fractal continuation is
applied to a fractal function $f$, a set of continuations is obtained. We
prove that, in the generic situation with $M=N=2$, this set of continuations
depends only on the function $f$ and is independent of the particular IFS $%
\mathcal{W}$ that was used to produce $f$. The proof relies on the detailed
geometrical structure of analytic fractal functions and on the Weierstass
preparation theorem.

\begin{figure}[htb]
\centering
\includegraphics[
natheight=5.694800in,
natwidth=13.416700in,
height=1.5342in,
width=4.5455in
]
{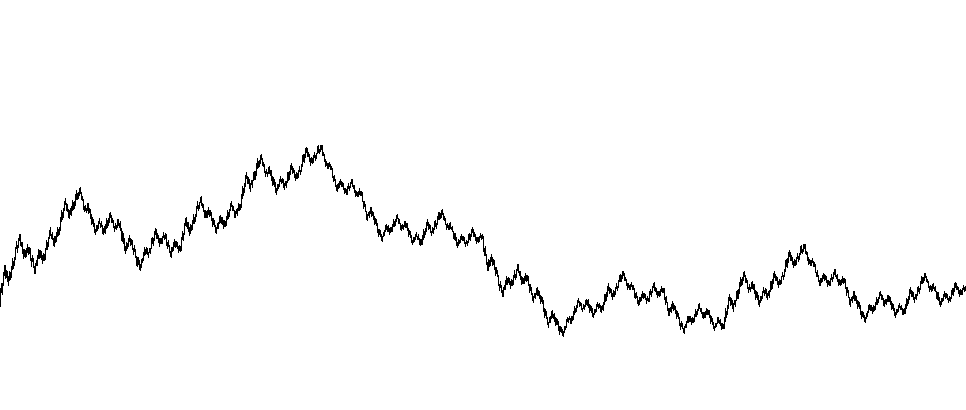}
\caption{This paper concerns analytic continuation, not only of analytic
functions, but also of non-differentiable functions such as the one whose
graph is illustrated here.}
\label{fif}
\end{figure}

The spirit of this paper is summarized in Figure \ref{fif}. Basic
terminology and background results related to iterated function systems
appear in Section~\ref{IFSsec}. In Section \ref{fifsec} we establish the
existence of \textit{fractal functions} whose graphs are the attractors of a
general class of IFS, which we call an \textit{interpolation IFS}. An 
\textit{analytic fractal function} is a fractal function whose graph is the
attractor of an analytic interpolation IFS. This includes the popular case
of affine fractal interpolation functions \cite{finterp}. An analytic
function is a special case of an analytic fractal function, as proved in
Section~\ref{analsec}. Fractal continuation, the main topic of this paper,
is introduced in Section \ref{continuesec}. The fractal continuation of an
analytic function is the usual analytic continuation. In general, however, a
fractal function defined on a compact domain, has infinitely many
continuations, this set of continuations having a fascinating geometric
structure as demonstrated by the examples that are also contained in Section %
\ref{continuesec}. The graph of a given fractal function can be the
attractor of many distinct analytic IFSs. We conjecture that the set of
fractal continuations of a function whose graph is the attractor of an
analytic interpolation IFS is independent of the particular IFS. Some cases
of this uniqueness result are proved in Section \ref{unicitysec}.

\section{Iterated Function Systems}

\label{IFSsec}

An \textbf{iterated function system (IFS)} 
\begin{equation*}
\mathcal{W=\{}\mathbb{X};w_{n},n\in\mathcal{I}\}
\end{equation*}
consists of a complete metric space $\mathbb{X\subset R}^{M}$ with metric $%
d_{\mathbb{X}},$ and $N$ continuous functions $w_{n}:\mathbb{X}\rightarrow 
\mathbb{X}$. The IFS $\mathcal{W}$ is called \textbf{contractive} if each
function $w$ in $\mathcal{W}$ is a contraction, i.e., if there is a constant 
$s\in\lbrack0,1)$ such that 
\begin{equation*}
d(w(x),w(y))\leq s \, d(x,y)
\end{equation*}
for all $x, y \in\mathbb{X}$. The IFS $\mathcal{W}$ is called an \textbf{%
invertible} IFS if each function in $\mathcal{W}$ is a homeomorphism of ${%
\mathbb{X}}$ onto ${\mathbb{X}}$. The definition of \textbf{analytic IFS} is
as given in the introduction. The IFS $\mathcal{W}$ is called an \textbf{%
affine }IFS if $\mathbb{X=R}^{M}$ and $w_{n}:\mathbb{R}^{M}\rightarrow%
\mathbb{R}^{M}$ is an invertible affine map for all $n\in\mathcal{I}$.
Clearly an affine IFS is analytic, and an analytic IFS is invertible.

The set of nonempty compact subsets of $\mathbb{X}$ is denoted $\mathbb{H} =%
\mathbb{H}(\mathbb{X})$. It is well-known that $\mathbb{H}$ is complete with
respect to the \textbf{Hausdorff metric} $h$, defined for all $S,T\in 
\mathbb{H}$, by 
\begin{equation*}
h(S,T)=\max \, \left \{\underset{s\in S}{\max} \, \underset{t\in T}{\min} \,
d_{\mathbb{X} }(s,t) ,\underset{t\in T}{\max}\, \underset{s\in S}{\min} \,
d_{\mathbb{X} }(s,t) \right\}\text{.}
\end{equation*}

\noindent Define $\mathcal{W} :\mathbb{H\rightarrow H}$ by 
\begin{equation}  \label{unioneq}
\mathcal{W} \mathcal{(}K)=\bigcup\limits_{n\in\mathcal{I}}w_{n}(K)
\end{equation}
for all $K\in\mathbb{H}$. Let $\mathcal{W}^{0}:\mathbb{H\rightarrow H}$ be
the identity map, and let $\mathcal{W}^{k}:\mathbb{H}\mathbb{\rightarrow H}$
be the $k$-fold composition of $\mathcal{W}$ with itself, for all integers $%
k>0.$

\begin{definition}
\label{attractordef} A set $A\in \mathbb{H}$ is said to be an \textbf{%
attractor} of $\mathcal{W}$ if ${\mathcal{W}} (A) = A$, and 
\begin{equation}  \label{attractorEq}
\lim_{k\rightarrow\infty}\mathcal{W}^{k}(K)=A
\end{equation}
for all $K\in\mathbb{H}$, where the convergence is with respect to the
Hausdorff metric.
\end{definition}

A basic result in the subject is the following \cite{hutchinson}.

\begin{theorem}
\label{Hutchthm} If IFS $\mathcal{W}$ is contractive, then $\mathcal{W}$ has
a unique attractor.
\end{theorem}

The remainder of this section provides the definition of a certain type of
IFS whose attractor is the graph of a function. We call this type of IFS an 
\textit{interpolation IFS}. We mainly follow the notation and ideas from 
\cite{finterp, FractalsEver, finterp2}.

Let $M\geq 2$ and $N$ an integer and $\mathcal{I=\{}1,2,...,N\}$. For a
sequence $x_{0}<x_{1}<\cdots <x_{N}$ of real numbers, let $L_{n}:\mathbb{%
R\rightarrow R}$ be the affine function and \linebreak $F_{n}:{\mathbb{R}}%
\times {\mathbb{R}}^{M-1}\rightarrow \ R^{M-1}, \; n=1,2,\dots ,N$ be a
continuous function satisfying the following properties:

\begin{enumerate}
\item[(a)] $L_{n}(x_0) = x_{n-1}$ and $L_n(x_N) = x_n$.

\item[(b)] There are points $y_{0}$ and $y_{N}$ in $\mathbb{R}^{M-1}$ such
that $F_{1}(x_{0},y_{0})=y_{0}$ and \linebreak $F_{N}(x_{N},y_{N})=y_{N}$.

\item[(c)] $F_{n+1}(x_0,y_0) = F_n(x_N,y_N)$ for $n= 1, 2, \dots, N-1$.
\end{enumerate}

Let $\mathcal{W}$ be the IFS 
\begin{equation}  \label{Weqn1}
\mathcal{W}=\{{\mathbb{R}}^M\, ; \, w_{n},\, n\in\mathcal{I\}}\text{, }
\end{equation}
where 
\begin{equation}  \label{Weqn2}
w_n (x,y) = (L_n(x), F_n(x,y)).
\end{equation}
Keeping condition (c) in mind, if we define, for each $n \in \mathcal{I}$, 
\begin{equation*}
y_n := F_{n+1}(x_0,y_0) = F_n(x_N,y_N),
\end{equation*}
then note that 
\begin{equation}  \label{interpeq}
w_{n}(x_{0},y_{0})=(x_{n-1},y_{n-1}) \qquad \text{and} \qquad
w_{n}(x_N,y_N)=(x_{n},y_{n}).
\end{equation}

\begin{definition}
\label{IntIFSdef} An \textbf{interpolation IFS} is an IFS of the form given
by (\ref{Weqn1}) and (\ref{Weqn2}) above that satisfies (1) $w_{n}$ is a
homeomorphism onto its image for all $n\in \mathcal{I}$, and (2) there is an 
$s\in \lbrack 0,1)$ and an $\mathcal{M}\in \lbrack 0,\infty )$ such that 
\begin{equation}
|F_{n}(x,y)-F_{n}(x^{\prime },y^{\prime })|\leq \mathcal{M}|x-x^{\prime
}|+s|y-y^{\prime }|  \label{conteq}
\end{equation}%
for all $x,x^{\prime }\in \mathbb{R},y,y^{\prime }\in {\mathbb{R}}^{M-1}$
and for all $n\in \mathcal{I}$. The term \textquotedblleft interpolation" is
justified by statement (2) in Theorem~\ref{functionthm} in the next section.
\end{definition}

\section{\label{fifsec} Fractal Functions and Interpolation}

Properties of an interpolation IFS are discussed in this section. Theorem~%
\ref{functionthm} is the main result.

\begin{lemma}
\label{contractionlem} If $\mathcal{W}$ is an interpolation IFS, then $%
\mathcal{W}$ is contractive with respect to a metric inducing the same
topology as the Euclidean metric on ${\mathbb{R}}^{M}$.
\end{lemma}

\begin{proof}
Let $L_{n}=a_{n}x+b_{n}$ and let $d$ be the metric on $\mathbb{R}^{M}$
defined by 
\begin{equation*}
d((x,y),(x^{\prime },y^{\prime }))=e|x-x^{\prime }|+|y-y^{\prime }|
\end{equation*}%
where $e\in (\mathcal{M}/(1-a),\infty )$ and $a=\max \{a_{n}:n\in \mathcal{%
I\}}$. The metric $d$ is a version of the "taxi-cab" metric and is
well-known to induce the usual topology on $\mathbb{R}^{M}$. Moreover, $%
w_{n} $ is a contraction with respect to the metric $d$: for $%
(x,y),(x^{\prime },y^{\prime })\in {\mathbb{R}}\times \mathbb{R}^{M-1}$, 
\begin{align*}
d(w_{n}(x,y),w_{n}(x^{\prime },y^{\prime }))& =e|L(x)-L(x^{\prime
})|+|F_{n}(x,y)-F_{n}(x^{\prime },y^{\prime })| \\
& \leq ea|x-x^{\prime }|+\mathcal{M}|x-x^{\prime }|+s|y-y^{\prime }|\text{ }
\\
& =\left( ea+\mathcal{M}\right) |x-x^{\prime }|+s|y-y^{\prime }| \\
& =ce|x-x^{\prime }|+s|y-y^{\prime }| \\
& \leq \max \{c,s\}\,d((x,y),(x^{\prime },y^{\prime }))
\end{align*}%
where $c=a+\mathcal{M}/e$ is monotone strictly decreasing function of $e$
for $e>0$, so $c$ is strictly less than $a+\mathcal{M}/(\mathcal{M}/(1-a))=1$
for $e\in (\mathcal{M}/(1-a),\infty ).$
\end{proof}

Theorem \ref{functionthm} generalizes results such as \cite[p.186, Theorem
5.4]{massopust}. Let $I=[x_{0},x_{N}]$, $C(I)=\{f:I\rightarrow \mathbb{R}%
^{M-1}:f$ is continuous$\}$ and $C_{0}(I):=\{g\in
C(I):g(x_{0})=y_{0},g(x_{N})=y_{N}\}.$

\begin{theorem}
\label{functionthm} If $\mathcal{W}$ is an interpolation IFS, then

\begin{enumerate}
\item The IFS $\mathcal{W}$ has a unique attractor $G:=G(f)$ that is the
graph of a continuous function $f \, :\, I \rightarrow {\mathbb{R}}^{M-1}$.

\item The function $f$ interpolates the data points $(x_0,y_0), (x_1, y_1),
\dots, (x_N, y_N)$, i.e., $f(x_n) = y_n$ for all $n$.

\item If $W \, : \, C_0(I) \rightarrow C_0(I)$ is defined by $(Wg)(x) =
F_{n}(L_{n}^{-1}(x),g(L_{n}^{-1}(x)))$ for $x\in [x_{n-1},x_{n}]$, for $n\in%
\mathcal{I}$, and for all $g\in C_0(I)$, then $W$ has a unique fixed point $%
f $ and 
\begin{equation*}
f = \lim_{n\rightarrow \infty} W^k(f_0)
\end{equation*}
for any $f_0 \in C_0(I)$.
\end{enumerate}
\end{theorem}

\begin{proof}
It is readily checked that the mapping $\mathcal{W}$ of equation~(\ref%
{unioneq}) takes $C_{0}(I)$ into $C_{0}(I)$ and also that the mapping $W$ of
statement (3) in Theorem~\ref{functionthm} takes $C_{0}(I)$ into $C_{0}(I).$
Moreover, if $G(f_0)$ is the graph of the function $f_0 \in C_0(I)$, then ${%
\mathcal{W}}(G(f_0) )= W(f_0)$. This implies, by property~(\ref{attractorEq}%
) of the attractor, that the function $f$ in statement (1), assuming that it
exists, is the same as the function $f$ of statement (3), assuming that it
exists. Statement (3) is proved first.

(3): That the map $W$ is a contraction on $C_0(I)$ with respect to the sup
norm can be seen as follows. For all $g\in C_{0}(I)$, 
\begin{align*}
|(Wg_{1})(x)-(Wg_{2})(x)| & \leq\max_{n\in\mathcal{I}}|F_{n}(L_{n}
^{-1}(x),g_{1}(L_{n}^{-1}(x)))-F_{n}(L_{n}^{-1}(x),g_{2}(L_{n}^{-1}(x)))| \\
& \leq\max_{n\in\mathcal{I}} s\,|g_{1}(L_{n}^{-1}(x))-g_{2}(L_{n}^{-1}(x))|
\\
& \leq s \,\left\vert g_{1}(x)-g_{2}(x)\right\vert,
\end{align*}
where $0 \leq s< 1$ is the constant in condition (\ref{conteq}). Statement
(3) now follows from the Banach contraction mapping theorem.

(1): According to Lemma~\ref{contractionlem}, the IFS $\mathcal{W}$ is
contractive. By Theorem~\ref{Hutchthm}, $\mathcal{W}$ has a unique attractor 
$G$. Let $G_0$ denote the graph of some function $f_0$ in $C_0(I)$. Using
statement (3) there is a function $f\in C_0(I)$ such that $f =
\lim_{k\rightarrow \infty} W^k (f_0)$. By what was stated in the first
paragraph of this proof and by the property~(\ref{attractorEq}) in the
definition of attractor, we have $G = \lim_{k\rightarrow \infty} {\mathcal{W}%
}^k(G_0) = G(\lim_{k\rightarrow \infty} W^k(f_0)) = g(f)$.

(2): The attractor $G$ must include the points $(x_{0},y_{0})$ and $%
(x_{N},y_{N})$ because they are fixed points of $w_{1}$ and $w_{N}.$ Hence,
by Equation~(\ref{interpeq}) $G$ must contain $(x_{n},y_{n}) = w_n(x_N,y_N)$
for all $n$.
\end{proof}

\begin{remark}
\label{genRem} Theorem~\ref{functionthm} remains true if $F_{n}:\mathbb{X}%
\rightarrow {\mathbb{R}}^{M-1}$ and $w_{n}:\mathbb{X}\rightarrow \mathbb{X}$%
, for all $n\in \mathcal{I}$, where $\mathbb{X}\subseteq {\mathbb{R}}^{M}$
is a complete subspace of ${\mathbb{R}}^{M}$, and $\mathbb{X}$ contains the
line segment $[x_{0},x_{N}]$ (treated as a subset of $\mathbb{R}^{M}$).
This, for example, is the situation in Theorem~\ref{corollary} in the next
section.
\end{remark}

\begin{definition}
\label{fractdef} A function $f$ whose graph is the attractor of an
interpolation IFS will be called a \textbf{fractal function}. A function $f$
whose graph is the attractor of an analytic interpolation IFS will be called
an \textbf{analytic fractal function}. Note that, although a fractal
function usually has the properties associated with a fractal set, there are
smooth cases. See the examples that follow.
\end{definition}

\begin{example}
\label{x2ex} (\textbf{parabola}) The attractor of the affine IFS 
\begin{equation*}
\mathcal{W=}\{\mathbb{R}^{2};w_{1}(x,y)=(x/2,y/4),w_{2}
(x,y)=((x+1)/2,(2x+y+1)/4)\}
\end{equation*}
is the graph $G$ of $f:[0,1]\rightarrow\mathbb{R}$, $f(x)=x^{2}$. For each $%
k\in\mathbb{N}$ the set of points $\mathcal{W}^{k}(\{(0,0),(1,1)\})$ is
contained in $G$ and the sequence $\left\{ \mathcal{W}^{k}(\{(0,0),(1,1)\})%
\right\} _{k=0}^{\infty}$ converges to $G$ in the Hausdorff metric. Also, if 
$f_{0}:[0,1]$ is a piecewise affine function that interpolates the data $%
\{(0,0),(1,1)\}$ then $f_k :=W^{k}(f_{0})$ is a piecewise affine function
that interpolates the data $\mathcal{W} ^{k}(\{(0,0),(1,1)\})$ and $\left\{
f_{k}\right\} _{k=0}^{\infty}$ converges to $f$ in $(C([0,1], \,d_{\infty}).$
The change of coordinates $(x,y)\rightarrow(y,x)$ yields an affine IFS whose
attractor is the graph of $f:[0,1]\rightarrow\mathbb{R}$, $f(x)=\sqrt{x}.$
\end{example}

\begin{example}
\label{x2fractalex} (\textbf{arc of infinite length}) Let $d_{1} +d_{2}>1,$ $%
d_{1},d_{2}\in(0,1)$. The attractor of the affine IFS $\mathcal{W}=\{\mathbb{%
R}^{2};w_{1},w_{2}\},$ where

\begin{equation*}
w_{1} 
\begin{bmatrix}
x \\ 
y%
\end{bmatrix}
=%
\begin{bmatrix}
a & 0 \\ 
c & d_{1}%
\end{bmatrix}%
\begin{bmatrix}
x \\ 
y%
\end{bmatrix}%
, \qquad w_{2} 
\begin{bmatrix}
x \\ 
y%
\end{bmatrix}
=%
\begin{bmatrix}
(1-a) & 0 \\ 
-c & d_{2}%
\end{bmatrix}%
\begin{bmatrix}
x \\ 
y%
\end{bmatrix}
+%
\begin{bmatrix}
2a \\ 
2c%
\end{bmatrix}
,
\end{equation*}
is the graph of $f:[0,2]\rightarrow\mathbb{R}$ which interpolates the data $%
(0,0),$ $(a,c),$ $(1,0)$ and has Minkowski dimension $D>1$ where (see e.g. 
\cite[p.204, Theorem 5.32]{massopust})%
\begin{equation*}
(a)^{D-1}d_{1}+(1-a)^{D-1}d_{2}=1\text{.}
\end{equation*}
\end{example}

\begin{example}
\label{smoothex} (\textbf{once differentiable function}) The attractor of
the affine IFS $\mathcal{W}=\{\mathbb{R}^{2};w_{1},w_{2}\},$ where

\begin{equation*}
w_{1} 
\begin{bmatrix}
x \\ 
y%
\end{bmatrix}
=%
\begin{bmatrix}
\frac{1}{3} & 0 \\ 
\frac{1}{2} & \frac{2}{9}%
\end{bmatrix}%
\begin{bmatrix}
x \\ 
y%
\end{bmatrix}
, \qquad w_{2} 
\begin{bmatrix}
x \\ 
y%
\end{bmatrix}
=%
\begin{bmatrix}
\frac{2}{3} & 0 \\ 
-\frac{1}{2} & \frac{2}{9}%
\end{bmatrix}%
\begin{bmatrix}
x \\ 
y%
\end{bmatrix}
+%
\begin{bmatrix}
\frac{2}{3} \\ 
1%
\end{bmatrix}
,
\end{equation*}
is the graph of a once differentiable function $f:[0,2]\rightarrow\mathbb{R}$
which interpolates the data $(0,0),$ $(2/3,1),$ $(2,0)$. The derivative is
not continuous. The technique for proving that the attractor of this IFS is
differentiable is described in \cite{berger}.
\end{example}

\begin{example}
\label{kigamiex} (\textbf{once continuously differentiable function}) The
attractor of the affine IFS $\mathcal{W}=\{\mathbb{R}^{2};w_{1},w_{2}\},$
where

\begin{equation*}
w_{1}%
\begin{bmatrix}
x \\ 
y%
\end{bmatrix}%
=%
\begin{bmatrix}
\frac{2}{5} & \frac{1}{5} \\ 
\frac{1}{5} & \frac{2}{5}%
\end{bmatrix}%
\begin{bmatrix}
x \\ 
y%
\end{bmatrix}%
,\qquad w_{2}%
\begin{bmatrix}
x \\ 
y%
\end{bmatrix}%
=%
\begin{bmatrix}
\frac{3}{5} & 0 \\ 
-\frac{1}{5} & \frac{1}{5}%
\end{bmatrix}%
\begin{bmatrix}
x \\ 
y%
\end{bmatrix}%
+%
\begin{bmatrix}
\frac{2}{5} \\ 
\frac{1}{5}%
\end{bmatrix}%
,
\end{equation*}%
is the graph of a continuous function $f:[0,1]\rightarrow \mathbb{R}$ that
possesses a continuous first derivative but is not twice differentiable, see 
\cite{barnfreiberg, berger}.
\end{example}

\begin{example}
\label{notdiffex}(\textbf{nowhere differentiable function of Weierstrass})
The attractor of the analytic IFS 
\begin{equation*}
\mathcal{W}=\{\mathbb{R}^{2};w_{1}(x,y)=(x/2,\xi y+\sin \pi
x),w_{2}(x,y)=((x+1)/2,\xi y-\sin \pi x)\}
\end{equation*}%
is the graph of $f:[0,1]\rightarrow \mathbb{R}$ well-defined, for $%
\left\vert \xi \right\vert <1$, by%
\begin{equation}
f(x)=\sum\limits_{k=0}^{\infty }\xi ^{k}\sin 2^{k+1}\pi x\text{ for all }%
x\in \lbrack 0,1],  \label{weiereqn}
\end{equation}%
This function $f(x)$ is not differentiable at any $x\in \lbrack 0,1]$, for
any $\xi \in \lbrack 0.5,1)$. See for example \cite[Ch. 5]{borwein}.
\end{example}

\section{Analytic Functions are Fractal Functions}

\label{analsec}

Given any analytic function $f:I\rightarrow\mathbb{R}$, we can find an
analytic interpolation IFS, defined on a neighbourhood $\mathbb{G}$ of the
graph $G(f)$ of $f$, whose attractor is $G(f)$. This is proved in two steps.
First we show that, if $f^{\prime}(x)$ is non zero and does not vary too
much over $I$, then a suitable IFS can be obtained explicitly. Then, with
the aid of an affine change of coordinates, we construct an IFS for the
general case.

\begin{lemma}
\label{analyticexthm} Let $f:[0,1]\rightarrow\mathbb{R}$ be analytic and
strictly monotone on $[0,1]$, with bounded derivative such that both (i) $%
\max_{x\in\lbrack0,1]}\left\vert f^{\prime}(x/2)/f^{\prime}(x)\right\vert <2$
and (ii) $\max_{x\in\lbrack0,1]}$ $\left\vert f^{\prime}((x+1)/2)\text{ }%
/f^{\prime}(x)\right\vert $ $<2$. Then there is a neighbourhood $\mathbb{%
G\subset R}^{2}$ of $G(f)$, complete with respect to the Euclidean metric on 
$\mathbb{R}^{2}$, such that 
\begin{equation*}
\mathcal{W=}\{\mathbb{G}%
;w_{1}(x,y)=(x/2,f(f^{-1}(y)/2)),w_{2}(x,y)=(x/2+1/2,f(f^{-1}(y)/2+1/2))\}
\end{equation*}
is an analytic interpolation IFS whose attractor is $G(f)$. 
\end{lemma}

\begin{proof}
First note that $G=G(f)$ is compact and nonempty. Also $G$ is invariant
under $\mathcal{W}$ because 
\begin{align*}
\mathcal{W}(G) & =f_{1}(G)\cup f_{2}(G) \\
& =\{(x/2,f(f^{-1}(f(x))/2):x\in\lbrack0,1]\}\cup \\
& \hskip 15mm \{(x/2+1/2,f(f^{-1}(f(x))/2+1/2):x\in\lbrack0,1]\} \\
& =\{(x/2,f(x/2)):x\in\lbrack0,1]\}\cup \{(x/2+1/2,f(x/2+1/2):x\in
\lbrack0,1]\} \\
& =G\text{.}
\end{align*}
Second, we show that property~(\ref{conteq}) holds on a closed neighbourhood 
$\mathbb{G}$ of $G.$ To do this we (a) show that it holds on $G$ and then
(b) invoke analytic continuation to get the result on a neighbourhood of $G$%
. Statement (a) follows from the chain rule for differentiation: for all $y$
such that $(x,y)\in G$ we have, 
\begin{align*}
\left\vert \frac{d}{dy}f(f^{-1}(y)/2)\right\vert & =\left\vert \frac{%
f^{^{\prime}}(f^{-1}(y)/2)}{2f^{^{\prime}}(f^{-1}(y))}\right\vert \\
& =\left\vert \frac{f^{\prime}(x/2)}{2 f^{\prime}(x)}\right\vert \quad (%
\text{since} \; y=f(x) \; \text{on} \; G) \\
& <1.
\end{align*}

To prove (b) we observe that both $v_{1}(y)=f(f^{-1}(y)/2)$ and $v_{2}
=f(f^{-1}(y)/2+1/2)$ are contractions for all $y$ such that $(x,y)\in G$ and
so, since they are both analytic, they are contractions for all $y$ in a
neighbourhood $\mathcal{N}$ of $f([0,1])$. Clearly $x/2$ and $x/2$ $+1/2$
are contractions for all $x\in\mathbb{R}$. Finally, it remains to show that $%
\mathbb{G\subset }\mathcal{N}$ can be chosen so that $w_{n}(\mathbb{%
G)\subset G}$ $(n=1,2)$. Let $\mathbb{G}$ be the union of all closed balls $%
B(x,y)$ of radius $\varepsilon>0$, centered on $(x,y)\in G,$ where $%
\varepsilon$ is chosen small enough that $\mathbb{G\subset}\mathcal{N}$.

The inverse of $w_{1}^{-1}:w_{1}(\mathbb{G)\rightarrow G}$ is well-defined
by $w_{1}^{-1}(x,y)=(2x,f^{-1}(2f(y)))$ and is continuous. Similarly we
establish that $w_{2}:\mathbb{G\rightarrow}\mathbb{G}$ is a homeomorphism
onto its image. Since $f$ and $f^{-1}$ are analytic, both $w_n$ and $%
w_n^{-1} $ are analytic, and hence $\mathcal{W}$ is analytic.
\end{proof}

\begin{theorem}
\label{corollary} If $f:I\rightarrow\mathbb{R}$ is analytic, then the graph $%
G(f)$ of $f$ is the attractor of an analytic interpolation IFS $\mathcal{W=\{%
} \mathbb{G};w_{1}(x,y),w_{2}(x,y)\}$, where $\mathbb{G}$ is a neighbourhood
of $G(f)$.
\end{theorem}

\begin{proof}
Let $f:I\rightarrow\mathbb{R}$ be analytic and have graph $G(f).$ We will
show, by explicit construction, that there exists an affine map $T:\mathbb{R}
^{2}\rightarrow\mathbb{R}^{2}$ of the form 
\begin{equation*}
T(x,y)=(ax+h,cx+dy)\text{, }ad\neq0,
\end{equation*}
such that $T(G(f))=G(g)$ is the graph of a function $g:[0,1]\rightarrow 
\mathbb{R}$ that obeys the conditions of Lemma \ref{analyticexthm} (wherein,
of course, you have to replace $f$ by $g$). Specifically, choose $L(x)=ax+h$
so that $T(I)=[0,1]$ and then choose the constants $c$ and $d$ so that 
\begin{equation*}
g(x)=f\left( \frac{x-h}{a}\right) d+\frac{c\left( x-h\right) }{a}
\end{equation*}
satisfies the conditions (i) and (ii) in Lemma \ref{analyticexthm}. To show
that this can always be done suppose, without loss of generality, that $%
I=[0,1]$ so that $a=1$ and $h=0,$ and let $d=1$. Then to satisfy condition
(i) and\ (ii) requires that 
\begin{equation*}
\max \left \{\max_{x\in\lbrack0,1]}\left\vert \frac{f^{\prime}\left(
x/2\right) +c}{f^{\prime}\left( x\right) +c}\right\vert ,\max_{x\in\lbrack
0,1]}\left\vert \frac{f^{\prime}\left( x/2+1/2\right) +c}{f^{\prime}\left(
x\right) +c}\right\vert \right \}<2,
\end{equation*}
which is true when we choose $c$ to be sufficiently large. Finally, let $%
\widetilde{\mathcal{W}}=\{\widetilde{\mathbb{G}};\widetilde{w}_{1} ,%
\widetilde{w}_{2}\}$ be the IFS, provided by Lemma \ref{analyticexthm}$,$
whose attractor is $G(g)$. Then $\mathcal{W}=\{\mathbb{G};w_{1}=T^{-1} \circ%
\widetilde{w_1}\circ T$, $w_{2}=T^{-1}\circ\widetilde{w}_{2}\circ T\},$
where $\mathbb{G}=T^{-1}(\widetilde{\mathbb{G}}),$ is an analytic
interpolation IFS whose attractor is the graph of $g$.

\end{proof}

The following example illustrates Theorem \ref{corollary}.

\begin{example}
\label{etothexex} (\textbf{the exponential function}) Consider the analytic
function $f(x) = e^x$ restricted to the domain $[1,2]$. An IFS obtained by
following the proof of Theorem~ \ref{corollary} is 
\begin{equation*}
\mathcal{W=}\left \{\mathbb{R}^{2};w_{1}(x,y)=(x/2+1/2,\sqrt{ey}
),w_{2}(x,y)=(x/2+1,e\sqrt{y}) \right \}.
\end{equation*}
Therefore the attractor of the analytic IFS $\mathcal{W}$ is the graph of $%
f:[1,2]\rightarrow\mathbb{R}$, $f(x)= e^x$. The change of coordinates $%
(x,y)\rightarrow(y,x)$ yields an analytic IFS whose attractor is an arc of
the graph of $\ln(x)$.
\end{example}

\section{\label{continuesec} Fractal Continuation}

This section describes a method for extending a fractal function beyond its
original domain of definiton. Theorem~\ref{mainthm} is the main result.

\begin{definition}
If $I\subset J$ are intervals on the real line and $f \, : \, I \rightarrow {%
\mathbb{R}}^{M-1}$ and $g \, : \, J \rightarrow {\mathbb{R}}^{M-1}$, then $g$
is called a \textbf{continuation} of $f$ if $f$ and $g$ agree on $I$.
\end{definition}

The following notion is useful for stating the results in this section. Let $%
\mathcal{I}^{\infty}$ denote the set of all strings $\sigma =\sigma_{1}
\sigma_{2}\sigma_{3} \cdots$, where $\sigma_k \in \mathcal{I}$ for all $k$.
The notation $\overline{\sigma_{1}\sigma_{2} \cdots\sigma_{m}}$ stands for
the periodic string $\sigma_{1}\sigma_{2} \cdots\sigma_{m}$ $\sigma_{1}
\cdots \sigma_{m}$ $\sigma_{1} \cdots$. For example, $\overline {12}=12121
\cdots$.

Given an IFS $\mathcal{W=}\{\mathbb{R}^{M};w_{n},n\in\mathcal{I}\}$, if $%
\sigma=\sigma_{1} \sigma_{2}\sigma_{3} \cdots \in\mathcal{I}^{\infty}$ and $%
k $ is a positive integer, then define $w_{\sigma|k}$ by 
\begin{equation*}
w_{\sigma|k}(x)=w_{\sigma _{1}}\circ w_{\sigma_{2}}\circ \cdots \circ
w_{\sigma_{k}}(x) =w_{\sigma_{1}}(w_{\sigma_{2}}( \dots (w_{\sigma_{k}}(x))
\cdots ))
\end{equation*}
for all $x\in {\mathbb{R}}^M$. Moreover, if each $w_n$ is invertible, define 
\begin{equation*}
w_{\theta|k}^{-1}:=w_{\theta_{1}}^{-1}\circ w_{\theta_{2}}^{-1}\circ...\circ
w_{\theta_{k}}^{-1}.
\end{equation*}
Note that, in general, $w_{\theta|k}^{-1} \neq (w_{\theta|k})^{-1}$.

A particular type of continuation, called a \textit{fractal continuation},
is defined as follows. Let $I$ be an interval on the real line and let $%
f:I\rightarrow\mathbb{R}^{M-1}$ be a fractal function as described in
Definition~\ref{fractdef}. In this section it is assumed that the IFS whose
attractor is $G(f)$ is invertible. Denote the inverse of $w_{n}(x,y)$ by 
\begin{equation*}
w_{n}^{-1}(x,y)=(L_{n}^{-1} (x),F_{n}^{\ast}(x,y)),
\end{equation*}
where $F_{n}^{\ast}:\mathbb{R}^{M}\rightarrow \mathbb{R}^{M-1}$ is, for each 
$n$, the unique solution to 
\begin{equation}  \label{eqInv}
F_{n}(L_{n}^{-1}(x),F_{n}^{\ast}(x,y))=y\text{.}
\end{equation}
Let $\theta \in \mathcal{I}^{\infty}$ and $G=G(f)$. Define 
\begin{equation}  \label{thetaeq1}
G_{\theta|k} :=w_{\theta|k}^{-1}(G)\text{.}
\end{equation}

\begin{prop}
\label{contProp} With notation as above, 
\begin{equation*}
G \subset G_{\theta|1}\subset G_{\theta|2}\subset \cdots \text{.}
\label{containeq}
\end{equation*}
Moreover, $G_{\theta|k}$ is the graph of a continuous function $f_{\theta|k}$
whose domain is 
\begin{equation*}
I_{\theta|k}:=L_{\theta|k}^{-1}(I) =L_{\theta_{1} }^{-1}\circ
L_{\theta_{2}}^{-1}\circ...\circ L_{\theta_{k}}^{-1}(I) .  \label{Ieqn}
\end{equation*}
\end{prop}

\begin{proof}
The inclusion $G_{\theta|k-1} \subset G_{\theta|k}$ is equivalent to $%
w_{\theta_k}(G) \subset G$, which follows from the fact that $G$ is the
attractor of $\mathcal{W}$. The second statement follows from the form of
the inverse as given in equation (\ref{eqInv}).
\end{proof}

It follows from Proposition~\ref{contProp} that 
\begin{equation}  \label{thetaeq2}
G_{\theta}:=\bigcup\limits_{k=0}^{\infty}G_{\theta|k}
\end{equation}
is the graph of a well-defined continuous function $f_{\theta}$ whose domain
is 
\begin{equation*}
I_{\theta}=\bigcup\limits_{k=0}^{\infty}I_{\theta|k}\text{ }.
\end{equation*}
Note that if $\theta \in \mathcal{I}^{\infty}$, then $f_{\theta}(x) =
f_{\theta|k}(x)$ for $x \in I_{\theta|k}$ and for any positive integer $k$.

\begin{definition}
The function $f_{\theta}$ will be referred to as the \textbf{fractal
continuation} of $f$ \textbf{with respect to} $\theta$, and $%
\{f_{_{\theta}}:\theta\in\mathcal{I}^{\infty}\}$ will be referred to as the 
\textbf{set of fractal continuations} of the fractal function $f$.
\end{definition}

Theorem \ref{mainthm} below states basic facts about the set of fractal
continuations. According to statement (3), the fractal continuation of an
analytic function is unique, not depending on the string $\theta \in {%
\mathcal{I}}^{\infty}$. Statement (4) is of practical value, as it implies
that stable methods, such as the chaos game algorithm, for computing
numerical approximations to attractors, may be used to compute fractal
continuations. The figures at the end of this section are computed in this
way.

\begin{theorem}
\label{mainthm} Let $\mathcal{W=}\{\mathbb{R}^{M};w_{n},n\in\mathcal{I}\}$
be an invertible interpolation IFS and let $G(f),$ the graph of $%
f:I\rightarrow\mathbb{R}^{M-1}$, be the attractor of $\mathcal{W}$ as
assured by Theorem \ref{functionthm}. If $\theta \in {\mathcal{I}}^{\infty}$%
, then the following statements hold:

\begin{enumerate}
\item 
\begin{equation*}
I_{\theta}=%
\begin{cases}
\mathbb{R} & \text{ if }\; \theta\in\mathcal{I}^{\infty} \setminus \{%
\overline {1},\overline{N}\}, \\ 
\mathbb{[}x_{0},\infty) & \text{ if }\; \theta=\overline{1}, \\ 
\mathbb{(-\infty},x_{N}] & \text{ if }\; \theta=\overline{N}\text{.}%
\end{cases}%
\end{equation*}

\item $f_{\theta}(x)=f(x)$ for all $x\in I$.

\item If $\mathcal{W}$ is an analytic IFS and $f$ is an analytic function on 
$I$, then $f_{\theta }(x)=\widetilde{f}(x)$ for all $x\in I_{\theta}$, where 
$\widetilde{f}:\mathbb{R\rightarrow R}^{M-1}$ is the (unique) real analytic
continuation of $f$.

\item For all $k\in\mathbb{N}$, the IFS 
\begin{equation*}
\mathcal{W}_{\theta|k}:=\{\mathbb{R}^{M};w_{\theta|k}^{-1}\circ
w_{n}\circ(w_{\theta|k}^{-1})^{-1}, \; n\in {\mathcal{I}} \}
\end{equation*}
has attractor $G_{\theta|k}=G(f_{\theta |k})$.
\end{enumerate}
\end{theorem}

\begin{proof}
(1): Each of the affine functions $L_n$ can be determined explicitely, and
it is easy to verify that $L_n^{-1}$ is an expansion. Moreover, the fixed
point of $L_n$ (and hence also of $L_n^{-1}$) lies properly between $x_0$
and $x_N$ for all $n$ except $n=1$ and $n=N$. The fixed point of $L_1$ is $%
x_0$ and the fixed point of $L_N$ is $x_N$. Statement (1) now follows from
the second part of Proposition~\ref{contProp}.

(2): It follows from Proposition~\ref{contProp} that $f_{\theta}(x)=f(x)$
for all $x\in I$, $\theta\in\mathcal{I}^{\infty}$.

(3): Since $F_n^*(x,y)$ and $L^{-1}_n$ are analyic for all $n$, each $%
w^{-1}_n$ is analytic. Therefore $f_{\theta}(x)$ is analytic and agrees with 
$\widetilde{f}(x)$ on $I.$ Hence $f_{\theta}(x)=\widetilde{f}(x)$, the
unique analytic continuation for $x\in I_{\theta}$.

(4): It is easy to check that condition (\ref{attractorEq}) in the
definition of attractor in Section~\ref{IFSsec} holds.
\end{proof}

\begin{cor}
\label{maincor} Let $\mathcal{W=}\{\mathbb{R}^{M};w_{n},n\in\mathcal{I}\}$
be an interpolaltion IFS, and let $G(f)$, the graph of $f:I\rightarrow%
\mathbb{R}^{M-1}$, be the attractor of $\mathcal{W}$ as assured by Theorem %
\ref{functionthm}. Let $f$ be analytic on $I$ and $\theta \in {\mathcal{I}}%
^{\infty}$.

\begin{enumerate}
\item If $\mathcal{W}$ is an affine IFS, then $f_{\theta}(x)=\widetilde{f}%
(x) $, where $\widetilde{f}:\mathbb{R\rightarrow R}^{M-1}$ is the real
analytic continuation of $f$.

\item If $M=2$ and $\mathcal{W}$ is the IFS constructed in Theorem~\ref%
{corollary}, then $f_{\theta}(x)=\widetilde{f}(x)$.
\end{enumerate}
\end{cor}

\begin{proof}
For both statements, ${\mathcal{W}}$ is analytic and hence the hypotheses of
statement (3) of Theorem~\ref{mainthm} are satisfied.
\end{proof}

\begin{remark}
This is a continuation of Remark~\ref{genRem}, and concerns a generalization
of Theorem~\ref{mainthm} to the case of an IFS in which the domain $\mathbb{X%
}$ of each $w_n : \mathbb{X }\rightarrow \mathbb{X}$ is a complete subspace
of ${\mathbb{R}}^{M}$. The attractor $G$ may not lie in the range of $%
w_{\theta|k}$ and the set-valued inverse $w_{\theta|k}^{-1}$ may map some
points and sets to the empty set. Nonetheless, it is readily established
that 
\begin{equation*}
G \subseteq G_{\theta|1}\subseteq G_{\theta|2}\subseteq \cdots
\end{equation*}
is an increasing sequence of compact sets contained in $\mathbb{X}$. We can
therefore define : 
\begin{equation*}
G_{\theta|k}=w_{\theta|k}^{-1}(G)\qquad \text{ and } \qquad
G_{\theta}=\bigcup \limits_{k=0}^{\infty}G_{\theta|k}\subset\mathbb{X}\text{.%
}
\end{equation*}
exactly as in Equations (\ref{thetaeq1}) and (\ref{thetaeq2}) and the
continuation $f_{\theta}$ as the function whose graph is $G_{\theta}$.
Theorem~\ref{mainthm} then holds in this setting.
\end{remark}

\begin{figure}[htb]
\fbox{\includegraphics[width=4.5in, keepaspectratio]{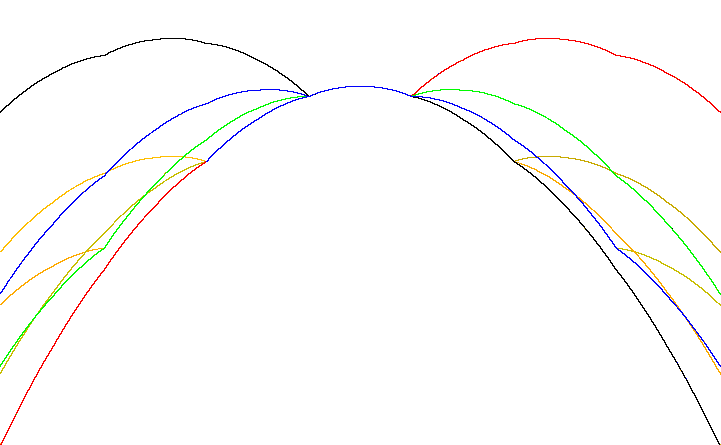}}
\caption{See Example \protect\ref{EXex}.}
\label{rainbow3}
\end{figure}

\begin{figure}[htb]
\vskip 4mm \includegraphics[width=2.5in, keepaspectratio]{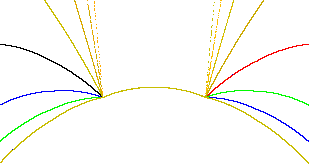} %
\vskip 3mm
\caption{Added detail for part of Figure \protect\ref{rainbow3} showing
additional continuations near the ends of the original function.}
\label{rainbow7}
\end{figure}

\begin{example}
\label{EXex} This example is related to Examples \ref{x2ex} and \ref%
{x2fractalex}. Let $G_{p}$ be the attractor of the affine IFS 
\begin{equation*}
\mathcal{W}_{p}=(\mathbb{R}^{2},w_{1}(x,y)=(0.5x,0.5x+py),w_{2}
(x,y)=(0.5x+1,-0.5x+py+1),
\end{equation*}
where $p\in(-1,0)\cup(0,1)$ is a parameter.

When $p=0.25$ the attractor $G_{0.25}$ is the graph of the analytic function 
$f\,:[0,2]\rightarrow {\mathbb{R}}$, where 
\begin{equation*}
f(x)=x(2-x)
\end{equation*}%
and, according to statement (3) of Theorem~\ref{mainthm} the unique
continuation is $f_{\theta }(x)=x(2-x)$ with domain $[0,\infty )$ if $\theta
=\overline{1}$, domain $(-\infty ,2]$ if $\theta =\overline{2}$, and domain $%
(-\infty ,\infty )$ otherwise.

\begin{figure}[htb]
\includegraphics[width=3.75in, keepaspectratio]{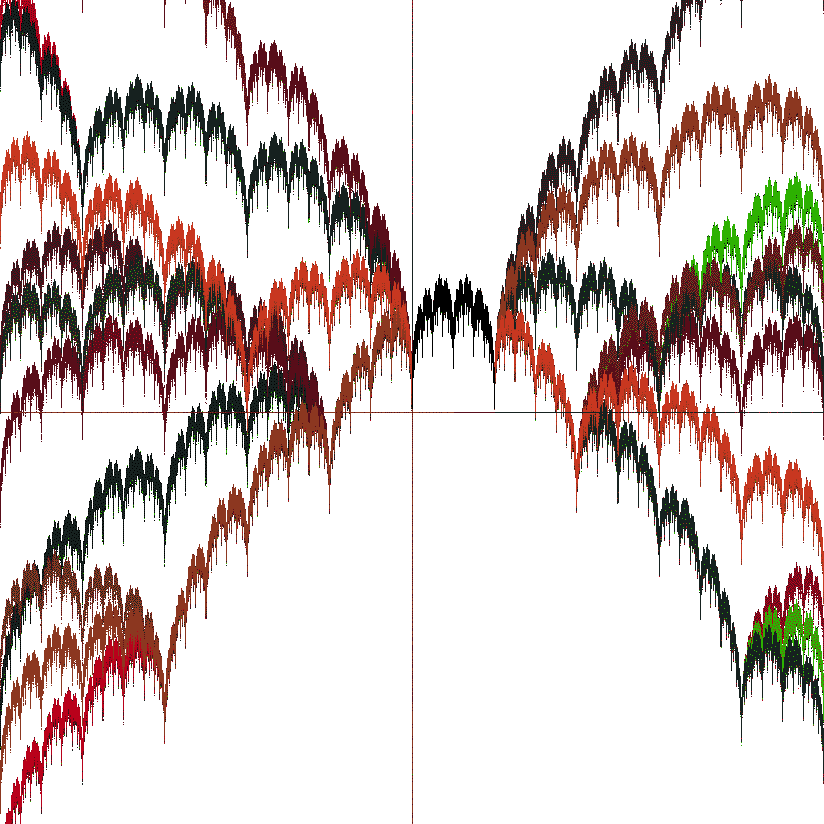}
\caption{See Example \protect\ref{EXex}. The colors help to distinguish the
different continuations of the fractal function whose graph is illustrated
in black near the center of the image. }
\label{iparabpt7}
\end{figure}

When $p=0.3$ the attractor is the graph of a non-differentiable function and
there are non-denumerably many distinct continuations $f_{\theta} \, : \,
(-\infty,\infty)\rightarrow\mathbb{R}$. Figure \ref{rainbow3} shows some of
these continuations, restricted to the domain $\lbrack-20,20]$. More
precisely, Figure \ref{rainbow3} shows the graphs of $f_{\theta|4}(x)$ for
all $\theta\in\{1,2\}^{\infty}.$ The continuation $f_{\overline{1}}(x)$, on
the right in black, coincides exactly, for $x\in\lbrack2,4]$, with \textit{%
all} continuations of the form $f_{1\rho}(x)$ with $\rho\in\{1,2\}^{\infty}.$
To the right of center: the blue curve is $G_{2111}$, the green curve is $%
G_{2211}$, and the red curve is $G_{2221}$. On the left: the lowest curve
(part red, part blue) is $G_{2222}$, the green curve is $G_{1222}$, the blue
curve is $G_{1122}$, and the black curve is $G_{1112}$. Also see Figure \ref%
{rainbow7}.

For $p=0.8$ the attractor $G_{0.8}$ is the graph of a fractal function $%
f_{0.8}$ whose graph has Minkowski dimension $(2-\ln(5/4)/\ln2)$. This graph 
$G_{0.8}$ is illustrated in the middle of Figure \ref{iparabpt7}. The window
for Figure \ref{iparabpt7} is $[-10,10]\times\lbrack-10,10]\subset \mathbb{R}%
^{2},$ and $f_{0.8}$ is the (unique) black object whose domain is $[0,2]$.
Figure \ref{iparabpt7} shows all continuations $f_{\theta|4}(x)$ for $\theta
\in \{1,2\}^{\infty}.$
\end{example}

\begin{figure}[tbh]
\centering
\fbox{%
\includegraphics[
natheight=14.221800in,
natwidth=14.221800in,
height=3.039in,
width=3.039in
]
{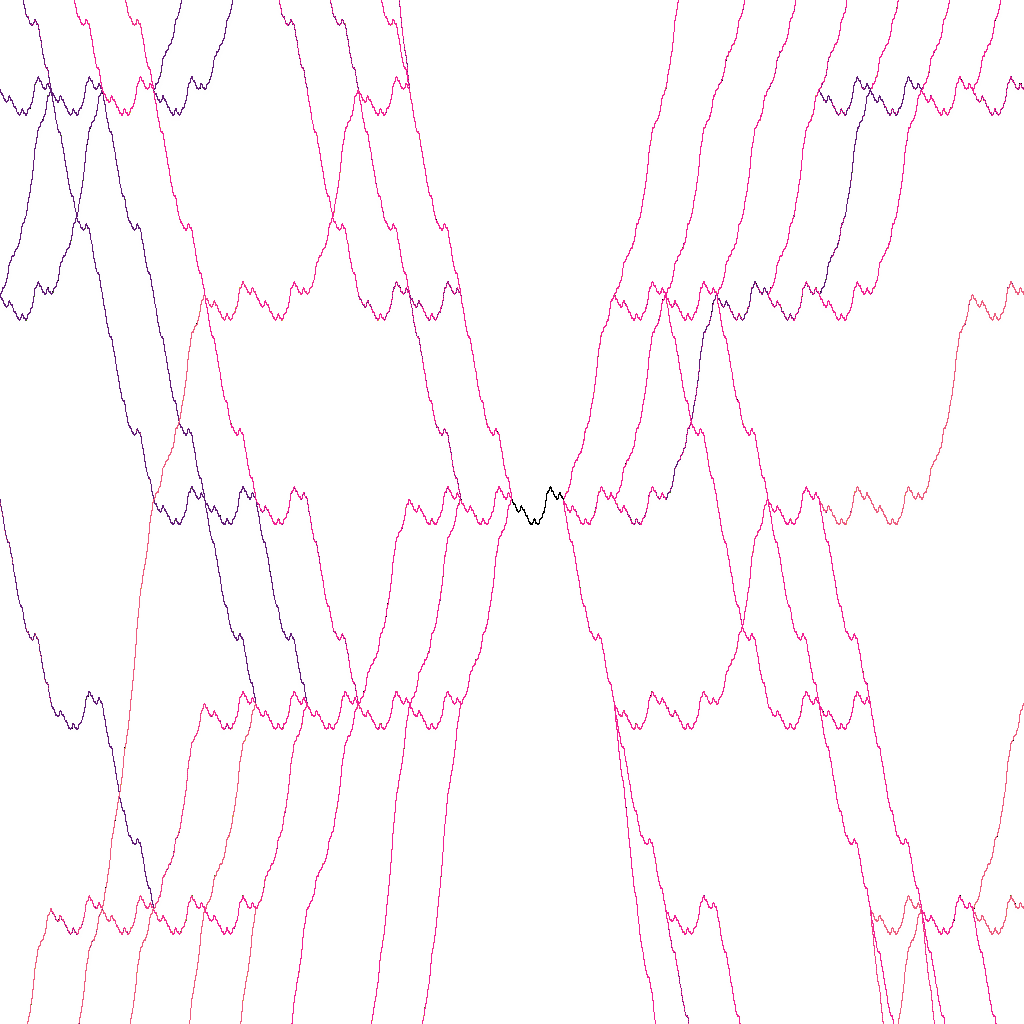}}
\caption{See Example \protect\ref{Ex9}.}
\label{ex9}
\end{figure}

\begin{example}
\label{notdiffexctd} (\textbf{continuation of a nowhere differentiable
function of Weierstrasse}) We continue Example \ref{notdiffex}. It is
readily calculated that, for all $\xi \in \lbrack 0,1)$, 
\begin{equation*}
w_{1}^{-1}(x,y)=(2x,\left( y-\sin 2\pi x\right) /\xi
),w_{2}^{-1}(x,y)=(2x-1,(y-\sin 2\pi x)/\xi )
\end{equation*}%
from which it follows that%
\begin{equation*}
f_{\theta }(x)=\sum\limits_{k=0}^{\infty }\xi ^{k+1}\sin 2^{k+1}\pi x
\end{equation*}%
with domain $[0,\infty )$ if $\theta =\overline{1}$, domain $(-\infty ,2]$
if $\theta =\overline{2}$, and domain $(-\infty ,\infty )$ otherwise. In
this example, all continuations agree, where they are defined, both with
each other and with unique function defined by periodic extension of
equation (\ref{weiereqn}).
\end{example}

When $\mathcal{W}$ in Theorem \ref{mainthm} is affine and $M=2$ write, for $%
n\in\mathcal{I}$, 
\begin{equation*}
F_{n}(x,y)=c_{n}x+d_{n}y+e_{n}.
\end{equation*}
We refer to the free parameter $d_{n}\in\mathbb{R}$, constrained by $%
\left\vert d_{n}\right\vert <1,$ as a \textbf{vertical scaling factor}. If
the verticle scaling factors are fixed and we require that the attractor
interpolate the data $\{(x_{i},y_{i})\}_{i=0}^{N}$, then the affine
functions $F_n$ are completely determined.

\begin{example}
\label{Ex9} The IFS $\mathcal{W}$ comprises the four affine maps $%
w_{n}(x,y)= \linebreak (L_{n} (x),F_{n}(x,y))$ that define the fractal
interpolation function $f:[0,1]\rightarrow\mathbb{R}$ specified by the data 
\begin{equation*}
\{(0,0.25),(0.25,0),(0.5,-0.25),(0.75,0.5),(1,0.25)\},
\end{equation*}
with vertical scaling factor $0.25$ on all four maps. Figure~\ref{ex9}
illustrates the attractor, the graph of $f$, together with graphs of all
continuations $f_{ijkl}$ where $i,j,k,l\in\{1,2,3,4\}$. The window is $%
\lbrack-10,10]^{2}$.
\end{example}

\begin{figure}[tbh]
\centering
\fbox{%
\includegraphics[
natheight=14.221800in,
natwidth=14.221800in,
height=3.039in,
width=3.039in
]
{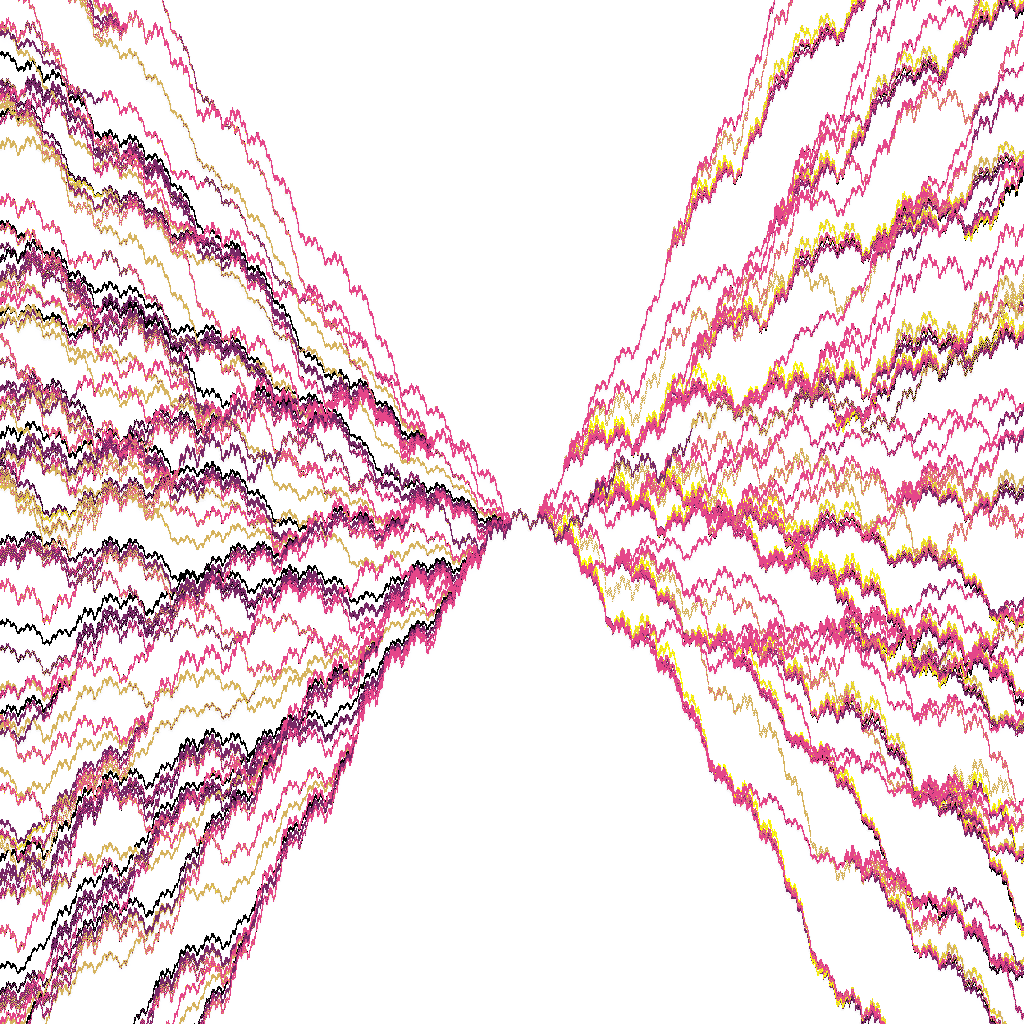} }
\caption{See Example~\protect\ref{paperfig2ex}.}
\label{paperfig2}
\end{figure}

\begin{example}
\label{paperfig2ex} The IFS $\mathcal{W}$ comprises the four affine maps
with respective vertical scaling factors $(0.55,0.45,0.45,0.45)$ such that
the attractor interpolates the data 
\begin{equation*}
\{(0,0.25),(0.25,0),(0.5,0.15),(0.75,0.6),(1,0.25)\}.
\end{equation*}
Figure~\ref{paperfig2} illustrates the attractor, an affine fractal
interpolation function $f:[0,1]\rightarrow\mathbb{R}$, together with all
continuations $f_{ijkl}$ where $i,j,k,l\in\{1,2,3,4\}$. The window is $%
\lbrack-20,20]^{2}$.
\end{example}

\begin{figure}[tbh]
\centering
\includegraphics[
natheight=14.221800in,
natwidth=14.221800in,
height=3.039in,
width=3.039in
]{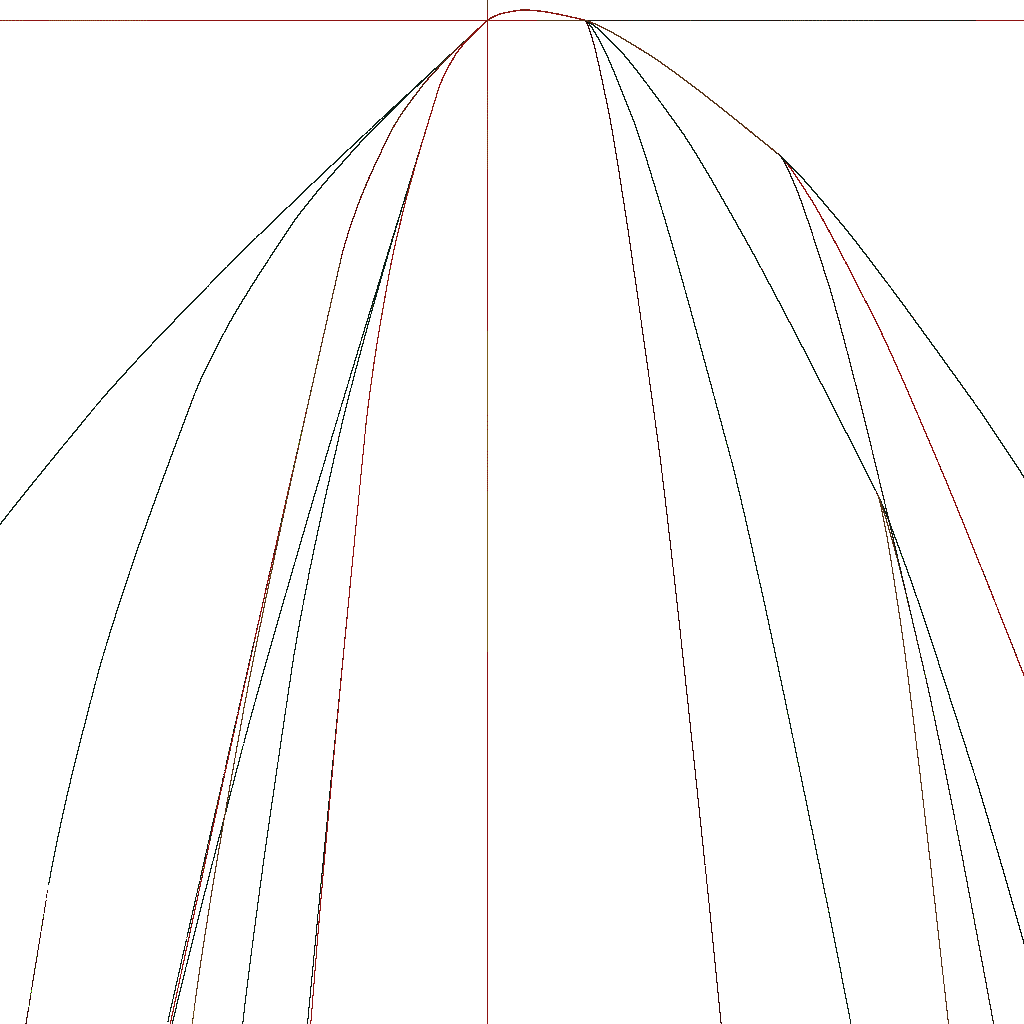}
\caption{See Example \protect\ref{smoothctdex}.}
\label{smooth01}
\end{figure}
\pagebreak

\begin{example}
\label{smoothctdex} Figure \ref{smooth01} shows continuations of the fractal
function $f:[0,2]\rightarrow \mathbb{R}$ described in Example \ref{smoothex}$%
.$ These are the continuations $f_{ijkl}(x)$ $(ijkl\in\{0,1\}^{4})$. The $x$%
-axis between $x=-10$ and $x=11$ and $y$-axis between $y=-100$ and $y=2$ are
also shown. The graph of $f$ is the part of the image above the $x$-axis$.$
\end{example}

\begin{figure}[hbt]
\centering
\includegraphics[
natheight=5.543400in,
natwidth=7.876700in,
height=2.4915in,
width=3.5276in
]
{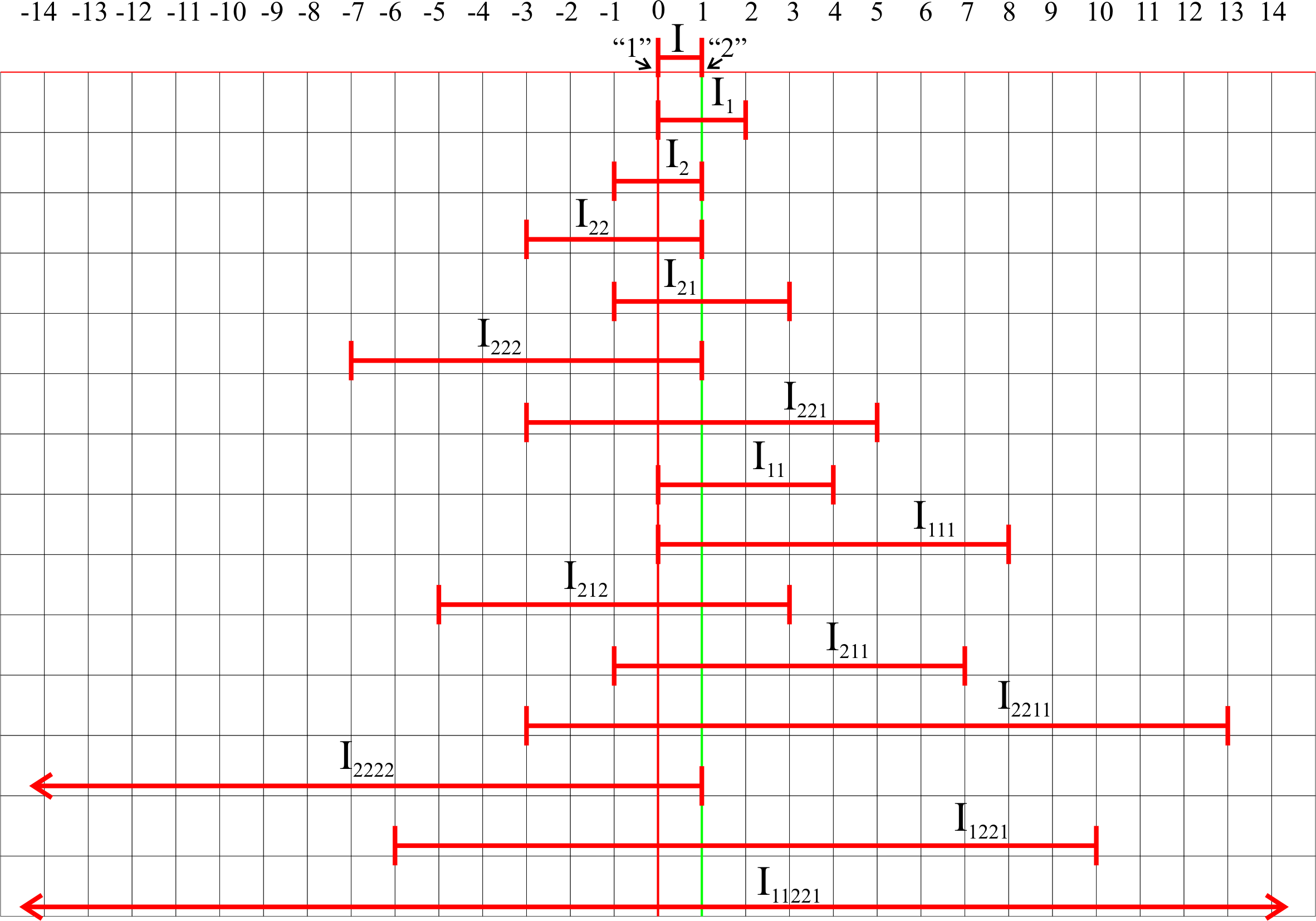}
\caption{An example of domains of agreements between ss analytic
continuations. See the end of Section \protect\ref{continuesec}}
\label{domains}
\end{figure}

In order to describes some relationships between the continuations $%
\{f_{\theta }:\theta \in \mathcal{I}^{\infty }\}$ (see the previous
examples), note that, for any finite string $\sigma $ and any $\theta ,
\theta^{\prime }\in \mathcal{I}^{\infty }$, 
\begin{equation*}
f_{\sigma \theta }(x)=f_{\sigma \theta ^{\prime }}(x)
\end{equation*}%
for all $x\in I_{\sigma }$. Consider the example $I=[0,1], \; N=2$, and 
\begin{equation*}
L_1(x) = \frac12 \, x, \qquad \qquad L_2(x) = \frac12 \, x + \frac12.
\end{equation*}
It is easy to determine $I_{\sigma }$ for various finite strings $\sigma $,
some of these intervals illustrated in Figure \ref{domains}. For example, we
must have 
\begin{equation*}
f_{22\theta }(x)=f_{22\theta ^{\prime }}(x)
\end{equation*}%
for all $x\in I_{22}=[-3,1]$ and for all $\theta ,\theta ^{\prime }\in 
\mathcal{I}^{\infty }$, but, as confirmed by examples, it can occur that $%
f_{212}(x)\neq f_{221}(x)$ for some $x\in I_{212}\cap I_{221}=[-3,3]$.

There is a natural probability measure on the collection of continuations on 
$\mathcal{I}^{\infty}$ defined by setting $\Pr(\theta_{i}=1)=0.5$ for all $%
i=1,2,..,$ independently. Then, because many continuations coincide over a
given interval, we can estimate probabilities for the values of the
continuations. For example, if $N=2$, $I=[0,1],$ and $a_{1}=a_{2}=1/2,$ then 
\begin{align*}
\Pr(f_{\theta}(x) & =f_{\overline{1}}(x) \, | \, x\in\lbrack1,2]) \geq 1/2;
\\
\Pr(f_{\theta}(x) & =f_{\overline{21}}(x) \, | \,x\in\lbrack1,2]) \geq 1/4;
\\
\Pr(f_{\theta}(x) & =f_{\overline{221}}(x) \, | \, x\in\lbrack1,2]) \geq1/8;
\\
\Pr(f_{\theta}(x) & =f\underset{n}{_{\underbrace{\overline{2...221}}} }(x)
\, | \, x\in\lbrack1,2]) \geq 1/2^{n}\text{ for all }n=1,2, \cdots
\end{align*}
In this sense, Figures \ref{rainbow3}, \ref{ex9}, and \ref{smooth01}
illustrate probable continuations.

\section{\label{unicitysec} Uniqueness of fractal continuations}

This section contains some results concerning the uniqueness of the set of
fractal continuations. Our conjecture is that an analytic fractal function $%
f $ has a unique set of continuations, indpependent of the particular IFS
that generates the graph of $f$ as the attractor. More precisely, suppose
that $G(f)$, the graph of a continuous function $f: I \rightarrow \mathbb{R}$%
, is the attractor of an analytic interpolation IFS $\mathcal{W}$ with set
of continuations $\{f_{\theta}:\theta\in\mathcal{I}^{\infty}\}$ as defined
in Section \ref{continuesec}, and the same $G(f)$ is also the attractor of
another analytic interpolation IFS $\widetilde{\mathcal{W} }=\{\mathbb{R}%
^{M};\widetilde{w}_{n},n\in \widetilde{\mathcal{I}}\}$ with set of
continuations $\{\widetilde{f}_{\theta}: \theta \in \widetilde{\mathcal{I}}%
^{\infty}\}$. The conjecture is that the two sets of continuations are equal
(although they may be indexed differently). This is clearly true if $f$ is
itself analytic, since an analytic function has a unique analytic
continuation. In this section we prove that the conjecture is true under
certain fairly general conditions when $f$ is not analytic. Recall that the
relevant IFSs are of the form 
\begin{equation}  \label{eq:form}
\begin{aligned} \mathcal{W} & =\left \{\mathbb{X\subset R}^{2}; \,
w_{n}(x,y)= (L_{n}(x),F_{n}(x,y)), \, n\in {\mathcal I} \right \}, \\ L_n(x)
& = a_n x + b_n \end{aligned}
\end{equation}
for $n=1,2, \dots ,N$. The first result concerns the extensions $f_{%
\overline{1}}(x)$ and $f_{\overline{N}}(x)$. This is a special case, but
introduces some key ideas.

\begin{theorem}
Let $\mathcal{W}$ and $\widetilde{\mathcal{W}}$ be analytic interpolation
IFSs, each with the same attractor $G(f)=G(\widetilde{f})$ but with possibly
different numbers, say $N$ and $\widetilde{N}$ of maps. Then 
\begin{equation*}
f_{\overline{1}}(x)=\widetilde{f}_{\overline{1}}(x) \qquad \text{ and }
\qquad f_{\overline{N}}(x)={\widetilde f}_{\overline{\widetilde{N}}}(x)
\end{equation*}
for all $x\in\mathbb{R}$ such that $(x,y)\in\mathbb{X}$ for some $y\in%
\mathbb{R}$.
\end{theorem}

\begin{proof}
As previously mentioned, it is sufficient to prove the theorem when $f$ is
not an analytic function. In this case $f$ must not possess a derivative of
some order at some point. By the self-similarity property~(\ref{eq:ss})
mentioned in the Introduction, $f$ must possess a dense set of such points.
Hence, as a consequence of the Weierstrass preparation theorem \cite%
{navascues}, if a real analytic function $g(x,y)$ vanishes on $G(f)$, then $%
g(x,y)$ must be identically zero. Now, since $L_{1}\circ \widetilde{L_{1}}=%
\widetilde{L_{1}}\circ L_{1}$ it follows, again from the self-similarity
property, that $\left( w_{1}\circ \widetilde{w_{1}}\right) (x,y)=\left( 
\widetilde{w_{1}}\circ w_{1}\right) (x,y)$ for all $(x,y)\in G(f).$ Then 
\begin{equation*}
g(x,y)=(w_{1}\circ \widetilde{w_{1}}-\widetilde{w_{1}}\circ w_{1})(x,y)
\end{equation*}%
vanishes on $G(f)$. Hence $w_{1}\circ \widetilde{w_{1}}=\widetilde{w_{1}}%
\circ w_{1}$ for all $(x,y)$ in $\mathbb{X}$. It follows, on multiplying on
the left by $w_{1}^{-1}$ and on the right by $w_{1}^{-1}$ that 
\begin{equation*}
w_{1}^{-1}\circ \widetilde{w_{1}}=\widetilde{w_{1}}\circ w_{1}^{-1},
\end{equation*}%
and similarly that $w_{1}\circ \widetilde{w_{1}}^{-1}=\widetilde{w_{1}}%
^{-1}\circ w_{1}$.

Now suppose that $(x,y)\in G(f_{\overline 1})$. Then $(x,y)\in G(f_{%
\overline {1|}k})=w_{1}^{-1}\circ...\circ w_{1}^{-1}(G(f))$ for some $k$.
Hence we can choose $l$ so large that 
\begin{align*}
\begin{array}{c}
\underbrace{\widetilde{w_{1}}\circ\widetilde{w_{1}}\circ \cdots \widetilde{%
w_{1}} } \\ 
l\text{ times}%
\end{array}
(x,y) & \in\widetilde{w_{1}}\circ\widetilde{w_{1}}\circ \cdots\widetilde{%
w_{1} }\circ w_{1}^{-1}\circ...\circ w_{1}^{-1}(G(f)) \\
& \subseteq w_{1}^{-1}\circ \cdots \circ w_{1}^{-1} 
\begin{array}{c}
\underbrace{\widetilde{w_{1}}\circ\widetilde{w_{1}}\circ \cdots \widetilde{%
w_{1}}} \\ 
l\text{ times}%
\end{array}
(G(f)) \\
& \subset G(f)
\end{align*}
which implies 
\begin{equation*}
(x,y)\in 
\begin{array}{c}
\underbrace{\widetilde{w_{1}}^{-1}\circ\widetilde{w_{1}}^{-1}\circ ....%
\widetilde{w_{1}}^{-1}} \\ 
l\text{ times}%
\end{array}
(G(f))
\end{equation*}
when $l$ is sufficiently large. Hence $G(f_{\overline 1})\subset G(%
\widetilde{f}_{\overline{\widetilde{1}}})$. The opposite inclusion is proved
similarly, as is the result for the other endpoint.
\end{proof}

\subsection{Differentiability of fractal functions.}

We are going to need the following result, which is interesting in its own
right, as it provides detailed information about analytic fractal functions.

\begin{theorem}
\label{lipschitzthm} Let $\mathcal{W}$ be an analytic interpolation IFS of
the form given in equation~(\ref{eq:form}) with attractor $G(f)$. Let $%
c\in\lbrack0,\infty)$ and $d\in\lbrack0,1)$ be real constants such that $%
\left\vert \frac{\partial }{\partial x}F_{n}(x,y)\right\vert <c$ and $0\leq
\left\vert \frac{\partial}{\partial y}F_{n}(x,y)\right\vert \leq d a_{n}$
for all $(x,y)$ in some neighborhood of $G(f)$, for all $n\in\mathcal{I}$.
The function $f:[x_{0} ,x_{N}]\rightarrow\mathbb{R}$ is lipschitz with
lipshitz constant $\lambda=a^{-1} c(1-d)^{-1},$ where $a=\min%
\{a_{n}:n=1,2,...,N\}$. That is: 
\begin{equation}
\left\vert f\left( s\right) -f(t)\right\vert \leq\lambda\, \left\vert
s-t\right\vert \text{ for all }s,t\in\lbrack x_{0},x_{N}].  \label{lipeqn}
\end{equation}
\end{theorem}

\begin{proof}
Consider the sequence of iterates $k=0,1,2, \dots$ 
\begin{equation*}
f_{k+1}(x)=\left( W f_{k}\right) (x)= F_{n}(L_{n}
^{-1}(x),f_{k}(L_{n}^{-1}(x)))
\end{equation*}
for $x \in [x_{n-1},x_{n}], n = 1,2, \dots, N$, where $W$ is as defined in
the statement of Theorem~\ref{functionthm}. Without loss of generality,
suppose that $\{f_{k}\}$ is contained in the neighborhood $\mathbb{X}$ of $%
G(f)$ mentioned in the statement of the theorem. It will first be shown, by
induction, that $f_k$ is lipschitz. Suppose that $f_{k}(x)$ is lipshitz on $%
[x_{0},x_{N}]$ with constant $\lambda $. Then, for all $s,t\in\lbrack
x_{n-1},x_{n}], \, n=1,2,...,N$ we have, by the self-replicating property,
by the mean value theorem for some $(\varsigma,\zeta )\in\mathbb{X}$, and by
the induction hypothesis that 
\begin{align*}
\left\vert f_{k+1}(s)-f_{k+1}(t)\right\vert & =\left\vert F_{n}(L_{n}
^{-1}(s),f_{k}(L_{n}^{-1}(s)))-F_{n}(L_{n}^{-1}(t),f_{k}(L_{n}^{-1}
(t)))\right\vert \\
& \leq\left\vert F_{n}(L_{n}^{-1}(s),f_{k}(L_{n}^{-1}(t)))-F_{n}(L_{n}
^{-1}(t),f_{k}(L_{n}^{-1}(t)))\right\vert \\
& +\left\vert F_{n}(L_{n}^{-1}(s),f_{k}(L_{n}^{-1}(s)))-F_{n}(L_{n}
^{-1}(s),f_{k}(L_{n}^{-1}(t)))\right\vert \\
& =\left\vert L_{n}^{-1}(s)-L_{n}^{-1}(t)\right\vert \cdot\left\vert \frac{%
\partial}{\partial x}F_{n}(x,y)\right\vert _{(\varsigma,t)} \\
& +\left\vert f_{k}(L_{n}^{-1}(s))-f_{k}(L_{n}^{-1}(t)))\right\vert
\cdot\left\vert \frac{\partial}{\partial y}F_{n}(x,y)\right\vert _{(s,\zeta)}
\\
& \leq\left( a_{n}^{-1} c+\lambda\, a_{n}^{-1}\, d\right) \left\vert
s-t\right\vert \leq\left( a^{-1}c+\lambda \, d\right) \left\vert
s-t\right\vert \\
& = \left( a^{-1}c+ a^{-1}c(1-d)^{-1} \, d\right) \left\vert s-t\right\vert
=\lambda\,\left\vert s-t\right\vert.
\end{align*}
Now suppose that $s<t$ and $s\in\lbrack x_{m-1},x_{m}]$ and $t\in\lbrack
x_{n-1},x_{n}]$ where $m<n$. Then 
\begin{align*}
\left\vert f_{k+1}(s)-f_{k+1}(t)\right\vert & \leq\left\vert f_{k+1}
(s)-f_{k+1}(x_{m})\right\vert +\left\vert f_{k+1}(x_{m})-f_{k+1}
(x_{m+1})\right\vert \\
& +\dots+\left\vert f_{k+1}(x_{n})-f_{k+1}(t)\right\vert \\
& \leq\lambda\,\left\vert s-x_{m}\right\vert +\lambda\,\left\vert
x_{m}-x_{m+1}\right\vert +\dots+\lambda\,\left\vert x_{n}-t\right\vert
=\lambda\, \left\vert s-t\right\vert.
\end{align*}
Therefore $f_{k}(x)$ is lipshitz on $[x_{0},x_{N}]$ with constant $\lambda$
for all $k$.

Now we use the fact that $\left\{ f_{k}\right\} $ converges uniformly to $f$
on $[x_{0},x_{N}]$, specifically 
\begin{equation*}
\begin{aligned} \max\{|f(x)-f_{k}(x)|:x\in\lbrack x_{0},x_{N}]\} &= \max\{
|W^kf(x)-W^kf_0(x)|: x\in\lbrack x_{0},x_{N}]\} \\ &\leq s^k \,
\max\{|f(x)-f_{0} (x)|:x\in\lbrack x_{0},x_{N}]\} \end{aligned}
\end{equation*}
for $s$ as in the proof of Theorem \ref{functionthm}. The unform limit of a
sequence of functions with lipschitz constant $\lambda$ is a lipshitz
function with constant $\lambda$.
\end{proof}

For an interpolation IFS $\mathcal{W}$ of the form given in equation~(\ref%
{eq:form}) with attractor $G(f)$, consider the IFS 
\begin{equation*}
\mathcal{L}:\mathcal{=\{}I=[x_{0},x_{N}];\,L_{n}(x)=a_{n}x+b_{n},n\in 
\mathcal{I}\},
\end{equation*}%
and let 
\begin{equation*}
\begin{aligned} \mathcal{D}_{\mathcal{L}}
&=\bigcup\limits_{k=0}^{\infty}\mathcal{L}
^{k}(\{x_{0},x_{1},\dots,x_{N}\})\backslash\{x_{0},x_{N}\}, \\ {\mathcal
D}_{\mathcal W} & = \{(x,f(x)) \, : \, x \in {\mathcal D}_{\mathcal L}\}.
\end{aligned}
\end{equation*}%
The set ${\mathcal{D}}_{\mathcal{W}}$ will be referred to as the set of 
\textbf{double points} of $G(f)$. The standard method for addressing the
points of the attractor of a contractive IFS \cite{FractalsEver} can be
applied to draw the following conclusions. If $(x,y)$ is a point of $G(f)$
that is not a double point, then there is a unique $\sigma \in {\mathcal{I}}%
^{\infty }$ such that 
\begin{equation*}
\begin{aligned} x & = \lim_{k\rightarrow \infty} L_{\sigma|k} (s), \\ (x,y)
& =\lim_{k\rightarrow \infty} w_{\sigma|k} (s,t), \end{aligned}
\end{equation*}%
where the limit is independent of $(s,t)\in \mathbb{X}$. If $(x,y)$ is a
double point, then there exist two distinct strings $\sigma $ such that the
above equations hold. In any case, we use the notation $\pi \,:\,{\mathcal{I}%
}^{\infty }\rightarrow {\mathbb{R}}$, 
\begin{equation*}
\pi (\sigma ):=\lim_{k\rightarrow \infty }L_{\sigma |k}(s),
\end{equation*}%
which is independent of $s\in \mathbb{R}$.

\begin{theorem}
\label{deriv1thm} Let $\mathcal{W}$ be an analytic interpolation IFS of the
form given in equation~(\ref{eq:form}) with attractor $G(f)$ and such that $%
0<\left\vert \partial F_{n}(x,y)/\partial y\right\vert <a_{n}$ for all $%
(x,y)\in\mathbb{X}$ and for all $n\in\mathcal{I}$. If $x$ is not a double
point of $G(f)$, then $f$ is differentiable at $x$. 
\end{theorem}

\begin{proof}
For now, fix $k\in\mathbb{N}.$ Define the function $Q:\mathbb{R}
^{k+1}\rightarrow\mathbb{R}$ by 
\begin{equation*}
Q(x_{1},x_{2}, \dots ,x_{k},y)=F_{\sigma_{1}}(x_{1},F_{\sigma_{2}}(x_{2} ,
\dots F_{\sigma_{k}}(x_{k},y) \dots)).
\end{equation*}
where the domain of $Q$ is all $(x_{1},x_{2}, \dots ,x_{k},y)$ such that $%
(x_{j},y)\in\mathbb{X}$ for all $j=1,2, \dots ,k.$ Let 
\begin{equation*}
H_{n}(x,y) := \frac{\partial}{\partial x}F_{n}(x,y), \qquad \qquad
K_{n}(x,y) := \frac{\partial}{\partial y}F_{n}(x,y)
\end{equation*}
for $n\in\mathcal{I}$. Because $F_n(x,y)$ is analytic for each $n$, there
are constants $c\in\lbrack0,\infty)$ and $d_{n} \in\lbrack0,a_{n})$ such
that 
\begin{equation}
\left\vert H_{n}(x,y)\right\vert \leq c, \qquad \qquad \left\vert
K_{n}(x,y)\right\vert <d_{n}  \label{bound2eqn}
\end{equation}
for all $n\in\mathcal{I}$ and for all $(x,y)$ in some neighborhood of $G(f)$%
. Using the notation 
\begin{align*}
Q_{x_{l}} (x_{1},x_{2},\dots,x_{k},y) & := \frac{\partial}{\partial x_{l}}%
Q(x_{1},x_{2},\dots,x_{k},y) \\
Q_{y}(x_{1} ,x_{2},\dots,x_{k},y) & :=\frac{\partial}{\partial y}%
Q(x_{1},x_{2},\dots,x_{k},y),
\end{align*}
we have the following for $l = 1,2, \dots ,k$: 
\begin{align*}
Q_{x_{l}}(x_{1},x_{2},\dots ,x_{k},y) &
=K_{\sigma_{1}}(x_{1},F_{\sigma_{2}}(x_{2},\dots
F_{\sigma_{k}}(x_{k},y)\dots))\cdot
K_{\sigma_{2}}(x_{2},F_{\sigma_{3}}(x_{3},\dots
F_{\sigma_{k}}(x_{k},y)\dots))\cdot\dots \\
& K_{\sigma_{l-1}}(x_{l-1},\dots F_{\sigma_{k}}(x_{k},y)\dots))\cdot
H_{\sigma _{l}}(x_{l},F_{\sigma_{l\mathfrak{+1}}}(x_{l+1},\dots
F_{\sigma_{k}}(x_{k},y)\dots));
\end{align*}
\begin{align*}
Q_{y}(x_{1},x_{2},\dots ,x_{k},y) &
=K_{\sigma_{1}}(x_{1},F_{\sigma_{2}}(x_{2},\dots
F_{\sigma_{k}}(x_{k},y)\dots))\cdot K_{\sigma_{2}}(x_{2},F_{\sigma
_{3}}(x_{3},\dots F_{\sigma_{k}}(x_{k},y)\dots))\cdot\dots \\
& K_{\sigma_{k-1}}(x_{k-1},F_{\sigma_{k}}(x_{k},y))\cdot K_{\sigma_{k}}
(x_{k},y).
\end{align*}
Using the intermediate value theorem repeatedly we have, for some $%
\eta_{l}\in\lbrack x_{l},x_{l}+\delta x_{l}]$, $l=1,2,\dots ,k$, and $\xi
\in\lbrack y,y+\delta y]$: 
\begin{align*}
\Delta Q := & \, Q(x_{1}+\delta x_{1},x_{2}+\delta x_{2},\dots ,x_{k}+\delta
x_{k},y+\delta y)-Q(x_{1},x_{2},\dots ,x_{k},y) \\
= & \, [Q(x_{1}+\delta x_{1},x_{2}+\delta x_{2},\dots ,x_{k}+\delta
x_{k},y+\delta y)-Q(x_{1},x_{2}+\delta x_{2},\dots ,x_{k}+\delta
x_{k},y+\delta y)]+ \\
& [Q(x_{1},x_{2}+\delta x_{2},\dots ,x_{k}+\delta x_{k},y+\delta y)-Q(x_{1}
,x_{2},x_{3}+\delta x_{3},\dots ,x_{k}+\delta x_{k},y+\delta y)]+\dots \\
& [Q(x_{1},x_{2},\dots x_{k-1},x_{k}+\delta x_{k},y+\delta y)-Q(x_{1}
,x_{2},\dots ,x_{k},y+\delta y)]+ \\
& [Q(x_{1},x_{2},\dots ,x_{k},y+\delta y)-Q(x_{1},x_{2},\dots ,x_{k},y)] \\
= & \, Q_{x_{1}}(\eta_{1},x_{2}+\delta x_{2},\dots ,x_{k}+\delta
x_{k},y+\delta y)\delta x_{1}+ \\
& Q_{x_{2}}(x_{1},\eta_{2},x_{3}+\delta x_{3}\dots ,x_{k}+\delta
x_{k},y+\delta y)\delta x_{2}+ \dots + \\
& Q_{x_{k}}(x_{1},x_{2},x_{3}\dots x_{k-1} ,\eta_{k},y+\delta y)\delta x_{k}+
\\
& Q_{y}(x_{1},x_{2}, x_{3}\dots , x_{k},\xi)\delta y\text{.}
\end{align*}
Let $\sigma = \pi^{-1}(x)$, which is well defined since $x$ is not a double
point, and \linebreak $x_{j}=\left( L_{\sigma|j} \right) ^{-1}(x).$ By the
self-repicating property~(\ref{unioneq}) 
\begin{equation*}
f_{k}(x)=F_{\sigma_{1}}(x_{1},F_{\sigma_{2}}(x_{2},\dots F_{\sigma_{k}}
(x_{k},f_{0}(x_k))\dots)).
\end{equation*}
Fix $k$ and $\sigma$ and let both $x$ and $(x+\delta x)$ lie in $L_{\sigma
|k}([0,1])$. Define 
\begin{equation*}
\begin{aligned} x_j & = \left( L_{\sigma|j} \right) ^{-1}(x) \\ y
&=f((L_{\sigma|k} )^{-1}(x)) \end{aligned} \qquad \begin{aligned}
x_{j}+\delta x_{j} &= (L_{\sigma|j} )^{-1}(x+\delta x)=(L_{\sigma|j}
)^{-1}(x)+(a_{\sigma_{1}}\dots a_{\sigma_{j} })^{-1}\delta x \\ y+\delta
y&=f((L_{\sigma|k} )^{-1}(x+\delta x)). \end{aligned}
\end{equation*}
Then 
\begin{align*}
& [f(x+\delta x)-f(x)]/\delta x \\
& =[F_{\sigma_{1}}(x_{1}+\delta x_{1},F_{\sigma_{2}}(x_{2}+\delta
x_{2},\dots F_{\sigma_{k}}(x_{k}+\delta x_{k},y+\delta y)\dots))- \\
& \hskip 5mm F_{\sigma_{1}}(x_{1},F_{\sigma_{2}}(x_{2},\dots
F_{\sigma_{k}}(x_{k},y)\dots))]/\delta x \\
& = \Delta Q / \delta x \\
& =[Q_{x_{1}}(\eta_{1},x_{2}+\delta x_{2},\dots ,x_{k}+\delta x_{k},y+\delta
y)\cdot(a_{\sigma_{1}})^{-1}+ \\
& \hskip 5mm Q_{x_{2}}(x_{1},\eta_{2},x_{3}+\delta x_{3}\dots ,x_{k}+\delta
x_{k},y+\delta y)\cdot(a_{\sigma_{1}}a_{\sigma_{2}})^{-1}+\dots \\
& \hskip 5mm Q_{x_{k}}(x_{1},x_{2},x_{3}\dots x_{k-1},\eta_{k},y+\delta
y)\cdot (a_{\sigma_{1}}\dots a_{\sigma_{k}})^{-1}]+ \\
& \hskip 5mm Q_{y}(x_{1},x_{2},x_{3}\dots
,x_{k},\xi)\cdot(a_{\sigma_{1}}\dots a_{\sigma_{k} })^{-1}\cdot \\
& \hskip 8mm \left [ f(\left( L_{\sigma|k}\right) ^{-1}(x)+((a_{\sigma_{1}
}\dots a_{\sigma_{k}})^{-1}\cdot \delta x)-f(\left( L_{\sigma|k}\right)
^{-1}(x)) \right ]/ ((a_{\sigma_{1}}\dots a_{\sigma_{k}})^{-1}\cdot \delta x)
\end{align*}
and 
\begin{equation*}
\begin{aligned} & [f(x+\delta x)-f(x)]/\delta x \\ &
=[c_{1}a_{\sigma_{1}}^{-1}+c_{2}d_{2,1}a_{\sigma_{1}}^{-1}a_{%
\sigma_{2}}^{-1}+
c_{3}d_{3,1}d_{3,2}a_{\sigma_{1}}^{-1}a_{\sigma_{2}}^{-1}a_{%
\sigma_{3}}^{-1}+\dots\nonumber\\ & \hskip 5mm +c_{k}d_{k,1}d_{k,2}\dots
d_{k,\left( k-1\right) }(a_{\sigma_{1}}\dots a_{\sigma_{k}})^{-1}] \\ &
=\frac{c_{1}}{a_{\sigma_{1}}}+\frac{c_{2}}{a_{\sigma_{1}}}\frac{d_{2,1}}{a_{%
\sigma_{2}}}+\dots\frac{c_{2}}{a_{\sigma_{1}}}\frac{d_{3,1}d_{3,2}}{a_{%
\sigma_{2}}a_{\sigma_{3}}}+\dots\frac{c_{k}}{a_{\sigma_{1}}}\frac
{d_{k,1}d_{k,2}\dots d_{k,\left( k-1\right)
}}{a_{\sigma_{2}}a_{\sigma_{3}}\dots a_{\sigma_{k}}}\\ & \hskip 5mm
+\frac{Q_{y}(x_{1},x_{2},x_{3}\dots ,x_{k},\xi)}{a_{\sigma_{1}}a_{\sigma_{2}
}a_{\sigma_{3}}\dots a_{\sigma_{k}}} \left [ \frac{f(
(L_{\sigma|k})^{-1}(x+\delta x))-f ( (L_{\sigma|k})^{-1}(x))}{
(L_{\sigma|k})^{-1}(x+\delta x)- (L_{\sigma|k})^{-1}(x)}\right ]
\end{aligned}
\end{equation*}
where 
\begin{align*}
c_{1} & =H_{\sigma_{1}}(\eta_{1},F_{\sigma_{2}}(x_{2}+\delta x_{2},\dots
F_{\sigma_{k}}(x_{k}+\delta x_{k},y+\delta y)\dots)), \\
c_{2} & =H_{\sigma_{2}}(\eta_{2},F_{\sigma_{3}}(x_{3}+\delta x_{3},\dots
F_{\sigma_{k}}(x_{k}+\delta x_{k},y+\delta y)\dots)), \\
& \ddots \\
c_{k} & =H_{\sigma_{k}}(\eta_{k},y+\delta y),
\end{align*}
and 
\begin{align*}
d_{2,1} &
=K_{\sigma_{1}}(x_{1},F_{\sigma_{2}}(\eta_{2},F_{\sigma_{3}}(x_{3}+\delta
x_{3},\dots F_{\sigma_{k}}(x_{k}+\delta x_{k},y)\dots))) \\
d_{3,1} & =K_{\sigma_{1}}(x_{1},F_{\sigma_{2}}(x_{2},F_{\sigma_{3}}(\eta
_{3},F_{\sigma_{3}}(x_{4}+\delta x_{4},\dots F_{\sigma_{k}}(x_{k}+\delta
x_{k},y)\dots)))) \\
d_{3,2} &
=K_{\sigma_{2}}(x_{2},F_{\sigma_{3}}(\eta_{3},F_{\sigma_{4}}(x_{4}+\delta
x_{4},\dots F_{\sigma_{k}}(x_{k}+\delta x_{k},y)\dots))) \\
& \ddots \\
d_{l,1} & =K_{\sigma_{1}}(x_{1},F_{\sigma_{2}}(x_{2},\dots
F_{\sigma_{l-1}}(x_{l-1},F_{\sigma_{l}}(\eta_{l},F_{\sigma_{l+1}}(x_{l+1}+%
\delta x_{l+1},\dots F_{\sigma_{k}}(x_{k}+\delta x_{k},y)\dots)))) \\
d_{l,2} & =K_{\sigma_{2}}(x_{2},F_{\sigma_{2}}(x_{2},\dots F_{\sigma_{l-1}
}(x_{l-1},F_{\sigma_{l}}(\eta_{l},F_{\sigma_{l+1}}(x_{l+1}+\delta
x_{l+1},\dots F_{\sigma_{k}}(x_{k}+\delta x_{k},y)\dots))))) \\
& \ddots \\
d_{l,l-1} & =F_{\sigma_{l-1}}(x_{l-1},F_{\sigma_{l}}(\eta_{l},F_{\sigma
_{l+1}}(x_{l+1}+\delta x_{l+1},\dots F_{\sigma_{k}}(x_{k}+\delta
x_{k},y)\dots))) \\
d_{k,1} & =K_{\sigma_{1}}(x_{1},F_{\sigma_{2}}(x_{2},\dots
F_{\sigma_{l-1}}(x_{l-1},F_{\sigma_{l}}(x_{l},F_{\sigma_{l+1}}(x_{l+1},\dots
F_{\sigma_{k}}(\eta_{k},y)\dots)))) \\
& \ddots \\
d_{k,k-1} & =K_{\sigma_{k-1}}(x_{k-1},F_{\sigma_{k}}(\eta_{k},y)).
\end{align*}
Note that the $c_{l}$ and $d_{l,m}$ depend explicitly on $k$ (which so far
is fixed). It follows that 
\begin{align*}
& \left\vert \frac{f(x+\delta x)-f(x)}{\delta x}-\left( \frac{c_{1} }{%
a_{\sigma_{1}}}+\frac{d_{2,1}}{a_{\sigma_{1}}}\frac{c_{2}}{a_{\sigma_{2}} }+ 
\frac{d_{3,1}d_{3,2}}{a_{\sigma_{1}}a_{\sigma_{2}}}\frac{c_{3}} {%
a_{\sigma_{3}}}+\cdots +\frac{d_{k,1}d_{k,2}\dots d_{k,\left( k-1\right) } }{%
a_{\sigma_{1}}a_{\sigma_{2}}\dots a_{\sigma_{k-1}}}\frac{c_{k}}{%
a_{\sigma_{k}} }\right) \right\vert \\
& \leq\left\vert \frac{Q_{y}(x_{1},x_{2},x_{3}\dots ,x_{k},\xi)}{%
a_{\sigma_{1} }a_{\sigma_{2}}a_{\sigma_{3}}\dots a_{\sigma_{k}}}\right\vert
\cdot\left\vert \frac{f(S_{k}(x+\delta x))-f(S_{k}(x))}{S^{k}(x+\delta
x)-S_{k} (x)}\right\vert \\
& \leq \lambda \, \left\vert \frac{Q_{y}(x_{1},x_{2},x_{3}\dots ,x_{k},\xi)}{%
a_{\sigma_{1} }a_{\sigma_{2}}a_{\sigma_{3}}\dots a_{\sigma_{k}}}\right\vert ,
\end{align*}
the last inequality by Theorem \ref{lipschitzthm}. The above is true for all 
$x,\left( x+\delta x\right) \in\lbrack x_{0},x_{N}], \linebreak \delta
x\neq0, k\in \mathbb{N}$. We also have 
\begin{equation*}
\begin{aligned} \left\vert \frac{Q_{y}(x_{1},x_{2},x_{3}\dots
,x_{k},\xi)}{a_{\sigma_{1} }a_{\sigma_{2}}a_{\sigma_{3}}\dots
a_{\sigma_{k}}}\right\vert & =\frac {\left\vert
K_{\sigma_{1}}(x_{1},F_{\sigma_{2}}(x_{2},\dots F_{\sigma_{k}}
(x_{k},\xi)\dots))\right\vert }{a_{\sigma_{1}}}\cdot \\ & \frac{\left\vert
K_{\sigma_{2}}(x_{2},F_{\sigma_{3}}(x_{3},\dots F_{\sigma_{k}}(x_{k}
,\xi)\dots))\right\vert }{a_{\sigma_{2}}}\cdots \\ & \frac{\left\vert
K_{\sigma_{k-1}}(x_{k-1},F_{\sigma_{k}}(x_{k} ,\xi))\right\vert
}{a_{\sigma_{k-1}}}\cdot\frac{\left\vert K_{\sigma_{k}
}(x_{k},\xi)\right\vert }{a_{\sigma_{k}}}\\ &
\leq\prod\limits_{j=1}^{k}\frac{d_{\sigma_{j}}}{a_{\sigma_{j}}}\leq
C^{k}\end{aligned}
\end{equation*}
for some $C \in\lbrack0,1)$, the last inequality by equation (\ref{bound2eqn}%
). Hence, for any $\varepsilon > 0$, we can choose $k$ so large that 
\begin{equation}  \label{eqEpsilon}
\left\vert \frac{f(x+\delta x)-f(x)}{\delta x}-\sum\limits_{m=1}^{k} \frac{%
c_{m}}{a_{\sigma_{m}}}\prod\limits_{l=1}^{m-1}\frac{d_{k,l}} {a_{\sigma_{l}}}%
\right\vert <\varepsilon/3 \text{.}
\end{equation}
Note that, by their definitions, for fixed $x$, the $c_{m}$s and $d_{k,l}$s
depend upon both $k$ and $\delta x$. Our next goal is to remove the
dependence on both $k$ and $\delta x$. For all $l$ and all $k\geq l$ define 
\begin{align}
C_{\sigma_{l}} &
:=H_{\sigma_{l}}(x_l,f(x_l)=H_{\sigma_{l}}(x_l,F_{\sigma_{l+1}}(x_{l+1},%
\dots F_{\sigma_{k}}(x_{k},f(x_k))\dots))  \label{comparisoneqn} \\
D_{\sigma_{l}} &
:=K_{\sigma_{l}}(x_l,f(x_l)=K_{\sigma_{l}}(x_l,F_{\sigma_{l+1}}(x_{l+1},%
\dots F_{\sigma_{k}}(x_{k},f(x_k))\dots)).  \notag
\end{align}
We are going to show that, for all $\varepsilon>0$ and for $\delta x$
sufficiently small, 
\begin{equation}
\left\vert \sum\limits_{m=1}^{k}\frac{C_{_{\sigma_{m}}}}{a_{\sigma_{m}}}%
\prod\limits_{l=1}^{m-1}\frac{D_{\sigma_{l}}}{a_{\sigma_{l}}}-\sum
\limits_{m=1}^{k}\frac{c_{m}}{a_{\sigma_{m}}}\prod\limits_{l=1}^{m-1}\frac{%
d_{k,l}}{a_{\sigma_{l}}}\right\vert <\varepsilon/3\text{,}  \label{eqnone}
\end{equation}
and that 
\begin{equation}
\left\vert \sum\limits_{m=1}^{\infty}\frac{C_{\sigma_{m}}}{a_{\sigma_{m}} }%
\prod\limits_{l=1}^{m-1}\frac{D_{\sigma_{l}}}{a_{\sigma_{l}}}-\sum
\limits_{m=1}^{k}\frac{C_{\sigma_{m}}}{a_{\sigma_{m}}}\prod\limits_{l=1}
^{m-1}\frac{D_{\sigma_{l}}}{a_{\sigma_{l}}}\right\vert <\varepsilon/3 \text{,%
}  \label{eqntwo}
\end{equation}
which taken together with inequality~(\ref{eqEpsilon}) imply 
\begin{equation}
\left\vert \frac{f(x+\delta x)-f(x)}{\delta x}-\sum\limits_{m=1}^{\infty}%
\frac{C_{\sigma_{m}}}{a_{\sigma_{m}}}\prod\limits_{l=1}^{m-1}\frac {%
D_{\sigma_{l}}}{a_{\sigma_{l}}}\right\vert < \varepsilon\text{.}
\label{threeeqn}
\end{equation}
For $|\delta x|$ sufficiently small 
\begin{align*}
& \left\vert \sum\limits_{m=1}^{k}\frac{C_{\sigma_{m}}}{a_{\sigma_{m}}}
\prod\limits_{l=1}^{m-1}\frac{D_{\sigma_{l}}}{a_{\sigma_{l}}}-\sum
\limits_{m=1}^{k}\frac{c_{m}}{a_{\sigma_{m}}}\prod\limits_{l=1}^{m-1} \frac{%
d_{k,l}}{a_{\sigma_{l}}}\right\vert \\
& \leq\sum\limits_{m=1}^{k}\frac{\left\vert C_{\sigma_{m}}-c_{m}\right\vert 
}{a_{\sigma_{m}}}\prod\limits_{l=1}^{m-1}\frac{\left\vert D_{\sigma_{l}
}\right\vert }{a_{\sigma_{l}}}+\sum\limits_{m=1}^{k}\frac{\left\vert
c_{m}\right\vert }{a_{\sigma_{m}}}\prod\limits_{l=1}^{m-1}\frac{\left\vert
D_{\sigma_{l}}-d_{k,l}\right\vert }{a_{\sigma_{l}}} < \varepsilon/3.
\end{align*}
The last inequality above follow from, for fixed $k,$ the continuous
dependence of the $c_{m}$s and $d_{k,l}$s on their independent variables,
and comparing $C_{\sigma_{m}}$ with $c_{m}$ and $D_{\sigma_{l}}$ with $%
d_{k,l}$ using the equalities (\ref{comparisoneqn}). (We need $|\delta x|$
small enough that $x+\delta x$ lies in $L_{\sigma|k}([0,1])$.) We have
established (\ref{eqnone}). Concerning inequality \linebreak (\ref{eqntwo}),
by equation~(\ref{bound2eqn}) the $c_{m}$'s are uniformly bounded and, for
some $(x,y)$, we have $\left\vert d_{k,l}\right\vert = |K_{\sigma_l}(x,y)|
\leq d_{\sigma_l} <a_{\sigma_l}$. Therefore $|d_{k,l}/a_{\sigma_l}| \leq
|d_{\sigma_l}/a_{\sigma_l}| \leq K$ for some constant $K < 1$. So inequality
(\ref{eqntwo}) follows from the absolute convergence of the series $%
\sum\limits_{m=1}^{\infty}\frac{C_{\sigma_{m}}}{a_{\sigma_{m}}}
\prod\limits_{l=1}^{m-1}\frac{D_{\sigma_{l}}}{a_{\sigma_{l}}}$. From
Equation (\ref{threeeqn}) it follows that 
\begin{equation*}
f^{\prime}(x) =\sum\limits_{m=1}^{k}\frac{C_{_{\sigma_{m}}}}{a_{\sigma _{m}}}%
\prod\limits_{l=1}^{m-1}\frac{D_{\sigma_{l}}}{a_{\sigma_{l}}}.
\end{equation*}
\end{proof}

Note that the last equality in the above proof actually provides a formula
for the derivative at each point that is not a double point.

\subsection{Unicity Theorem}

We conjecture that the uniqueness of the set of continuations holds in
general. The following theorem provides a proof in ${\mathbb{R}}^2$ under
the assumption that the derivative $f^{\prime }(x)$ does not exist at all
points $x$, although we conjecture that uniqueness holds in ${\mathbb{R}}^M,
M\geq 2$, and it is sufficient to assume that $f(x)$ is not analytic. It is
also assumed that there is a bound $\left\vert \partial F_{n}(x,y)/\partial
y\right\vert <a_{n}$, where the $a_n$ are as given in equation~(\ref{eq:form}%
). As an example, consider the case of affne fractal interpolation
functions, where $F_{n}(x,y)=(a_{n}x+b_{n},c_{n}x+d_{n}y+g_{n})$. Then for
Theorem \ref{unique1thm} to apply we need $\left\vert d_{n}\right\vert
<a_{n} $ for all $n$.

\begin{theorem}
\label{unique1thm} Let $\mathcal{W }= \{\mathbb{X\subset R}^{2}; \, w_{n}
(x,y)=(L_{n}(x),F_{n}(x,y)), \,n\in\mathcal{I}\}$ and \linebreak $\widetilde{%
\mathcal{W}}=\{\mathbb{X\subset R}^{2}; \, \widetilde{w}_{n}(x,y))=(%
\widetilde{L} _{n}(x),\widetilde{F}_{n}(x,y)), \, n\in\mathcal{I}\}$ be
analytic interpolation IFSs as in equation (\ref{eq:form}) such that $%
0<\left\vert \partial F_{n}(x,y)/\partial y\right\vert <a_{n}$ and $%
0<\left\vert \partial\widetilde{F}_{n}(x,y)/\partial y\right\vert <%
\widetilde{a}_{n}$ for all $(x,y)\in\mathbb{X}$, for all $n\in\mathcal{I}$.
If both $\mathcal{W}$ and $\widetilde{\mathcal{W}}$ have the same attractor $%
G(f)$ such that $f^{\prime }(x)$ does not exist at $x=x_{n},$ for all $%
n=0,1,2,\dots ,N$, then $\mathcal{W=}\widetilde{\mathcal{W}}$.
\end{theorem}

\begin{proof}
For simplicity we restrict the proof to the case $N=2$. The proof of the
result for arbitrary many interpolation points is similar.

We first prove that the set of double points of $G(f)$ with repect to $%
\widetilde{\mathcal{W}}$ is the same as the set of double points of $%
\mathcal{W}$. The interpolation points for $\mathcal{W}$ are $\{0,x_{1},1\}$
and the interpolation points for $\widetilde{\mathcal{W}}$ are $\{0,%
\widetilde{x}_{1},1\}$. By Theorem~\ref{deriv1thm} $f(x)$ is differentiable
at all points that are not double points with respect to $\mathcal{W}$ and
also at all points that are not double points with respect to $\widetilde{%
\mathcal{W}}$. Moreover, $f(x)$ is not differentiable at all double points
with respect to $\mathcal{W}$ and also not differentiable at all points
which are double points with respect to $\widetilde{\mathcal{W}}$.
(Otherwise $f(x)$ must be differentiable at $x_{1}$ which would imply that $%
f(x)$ is differentiable everywhere, contrary to the assumptions of the
theorem.) It follows that $f(x)$ is not differentiable at $x$ if and only if 
$x$ is a double points with respect to $\mathcal{W}$ if and only if $x$ is a
double point with respect to $\widetilde{\mathcal{W}}$.

We next prove that $w_{n}(x,y)=\widetilde{w}_{n}(x,y)$ for all $(x,y)\in
G(f) $ and $n=1,2,\dots ,N$. Since $\widetilde{x}_{1}$ is a double point of $%
G(f)$ with respect to $\mathcal{W}$ there must be $\sigma|k \neq\emptyset$
such that $w_{\sigma|k}(x_{1},f(x_{1}))=(\widetilde{x}_{1},f(\widetilde {x}%
_{1}))$. Since $x_{1}$ is a double point of $G(f)$ with respect to $%
\widetilde{\mathcal{W}}$ there must be $\widetilde \sigma|\widetilde k$ such
that ${\widetilde w}_{\widetilde {\sigma}|\widetilde{k}}(\widetilde{x}_{1},f(%
\widetilde{x}_{1}))=(x_{1} ,f(x_{1}))$. It follows that ${\widetilde w}_{(%
\widetilde{\sigma}|\widetilde{k} )}(w_{(%
\sigma|k)}(x_{1},f(x_{1})))=(x_{1},f(x_{1}).$ Since $\widetilde{w}_{(%
\widetilde {\sigma}|\widetilde{k})}\circ w_{(\sigma|k)}:G(f)\rightarrow G(f)$%
, we can write $\widetilde{w}_{(\widetilde{\sigma}|\widetilde{k} )}\circ
w_{(\sigma|k)}(x,y)=\overline{w}(x,y)=(\overline{L}(x),\overline {F}(x,y))$
where, similar in form to the functions $\widetilde{w}_{\widetilde{n} }(x,y)$
and $w_{n}(x,y)$ that comprise the two IFSs, $\overline{L} (x)=\overline{a}x+%
\overline{h}$ is a real affine contraction and $\overline {F}(x,y)$ is
analytic in a neighborhood of $G(f)$ and has the property, by the chain
rule, that $\left\vert \frac{\partial\overline{F}}{\partial y}(x,y)\}
\right\vert < \overline{a}$ in a neighborhood of $G(f)$. It is also the case
that $\overline{a} x_{1}+\overline{h}=x_{1}$ and $\overline{F}({\overline L}%
^{-1}(x),f({\overline L}^{-1} (x)))=f(x)$ in a neighborhood of $x_{1}$ and ${%
\overline L}(x_1) = x_1, \; \overline{F}(x_1,f(x_{1} ))=f(x_{1})$. Using the
analyticity of $F(x,y)$ in $x$ and $y$, 
\begin{align*}
\frac{f(x_{1}+\delta x)-f(x_{1})}{\delta x} & =\frac{\overline{F}
(L^{-1}(x_{1}+\delta x),f(L^{-1}(x_{1}+\delta x)))-\overline{F}(L^{-1}
(x_{1}),f(L^{-1}(x_{1})))}{\delta x} \\
& =\overline{F}_{x}(x_{1},f(x_{1}))\overline{a}^{-1} \\
& +\overline{F}_{y} (x_{1},f(x_{1}))\overline{a}^{-1}\frac{\left( f(x_{1}+%
\overline{a}^{-1}\delta x)-f(x_{1})\right) }{\overline{a}^{-1}\delta x}%
+o(\delta x).
\end{align*}
This implies that the following limit exists:%
\begin{align*}
& \lim_{\delta x\rightarrow0}\left\{ \frac{f(x_{1}+\delta x)-f(x_{1})}{%
\delta x}-\overline{F}_{y}(x_{1},f(x_{1}))\overline{a}^{-1}\frac{\left(
f(x_{1}+\overline{a}^{-1}\delta x)-f(x_{1})\right) }{\overline{a}^{-1}\delta
x}\right\} \\
& =(1-\overline{F}_{y}(x_{1},f(x_{1}))\overline{a}^{-1})f^{\prime}(x_{1})=%
\overline{F}_{x}(x_{1},f(x_{1}))\overline{a}^{-1},
\end{align*}
which implies 
\begin{equation*}
f^{\prime}(x_{1})=\frac{\overline{F}_{x}(x_{1},f(x_{1}))}{(\overline {a}-%
\overline{F}_{y}(x_{1},f(x_{1})))}.
\end{equation*}

We have shown that if $\sigma|k \neq\emptyset$ then $f(x)$ is differentiable
at $x_{1}$, which is not true. Therefore $\sigma|k =\emptyset$ which implies 
$x_{1}=\widetilde{x}_{1}$ and hence $w_{n}(x,y)=\widetilde {w}_{n}(x,y)$ for
a dense set of points $(x,y)$ on $G(f)$. It follows that $w_{n}(x,y)=%
\widetilde {w}_{n}(x,y)$ for all $(x,y)\in G(f)$ and $n=1,2.$

To show that $\mathcal{W=}\widetilde{\mathcal{W}}$, i.e., that $w_{n}(x,y)-%
\widetilde{w}_{n}(x,y)$ for all $(x,y)\in \mathbb{X}$, define an analytic
function of two variables, $a:\mathbb{X\rightarrow R}$ by $%
a(x,y):=w_{n}(x,y)-\widetilde{w}_{n}(x,y)$ for all $(x,y)\in \mathbb{X}$. It
was shown above that $a(x,y)=0$ for all $(x,y)\in G(f)$. That $a(x,y)=0$ for
all $(x,y)\in \mathbb{X}$ follows from the Weierstass preparation theorem 
\cite{narasimhan}.
\end{proof}

\begin{center}
\textbf{AKNOWLEDGEMENT}
\end{center}

We thank Louisa Barnsley for help with the illustrations.


\begin{thebibliography}{99}
\bibitem{finterp} M. F. Barnsley, Fractal functions and interpolation\textit{%
, Constr. Approx.} \textbf{2 }(1986) 303-329.

\bibitem{FractalsEver} M. F. Barnsley, \textit{Fractals Everywhere},
Academic Press, 1988; 2nd Edition, Morgan Kaufmann 1993; 3rd Edition, Dover
Publications, 2012.

\bibitem{finterp2} M. F. Barnsley and A. N. Harrington, The calculus of
fractal interpolation functions\textit{, Journal of Approximation Theory} 
\textbf{57} (1989) 14-34.

\bibitem{barnfreiberg} M. F. Barnsley, U. Freiberg, Fractal transformations
of harmonic functions\textit{, Proc. SPIE} \textbf{6417 }(2006).




\bibitem{berger} M.A. Berger, Random affine iterated function systems: curve
generation and wavelets, \textit{SIAM Review }\textbf{34} (1992) 361-385.

\bibitem{borwein} D.H. Bailey, J.M. Borwein, N.J. Calkin, R. Girgensohn,
D.R. Luke, V.H. Moll, \textit{Experimental Mathematics in Action, }A.K.
Peters, 2006.

\bibitem{hutchinson} J. E. Hutchinson, Fractals and self-similarity\textit{,
Indiana Univ. Math. J.} \textbf{30 }(1981) 713--747.


\bibitem{massopust0} Peter Massopust, \textit{Fractal Functions, Fractal
Surfaces, and Wavelets, }Academic Press, New York, 1995.

\bibitem{massopust} Peter Massopust, \textit{Interpolation and Approximation
with Splines and Fractals,} Oxford University Press, Oxford, New York, 2010.

\bibitem{narasimhan} R. Narasimhan, Introduction to the theory of analytic
spaces, \textit{Lecture Notes in Mathematics, }volume 25, Springer, 1966.

\bibitem{navascues} M. A. Navascues, Fractal polynomial interpolation\textit{%
, Zeitschrift f\"{u}r Analysis u. i. Anwend}, \textbf{24} (2005) 401-414.

\bibitem{prasad} Srijanani Anurag Prasad, \textit{Some Aspects of
Coalescence and Superfractal Interpolation, }Ph.D Thesis, Department of
Mathematics and Statistics, Indian Institute of Technology, Kanpur, March
2011.

\bibitem{scealy} Robert Scealy, $V$\textit{-variable fractals and
interpolation}, Ph.D. Thesis, Australian National University, 2008

\bibitem{tosan} Eric Tosan, Eric Guerin, Atilla Baskurt, Design and
reconstruction of fractal surfaces\textit{. }In IEEE Computer Society,
editor, \textit{6th International Conference on Information Visualisation IV
2002, London, UK }pp. 311-316, July 2002.

\bibitem{tricot} Claude Tricot, \textit{Curves and Fractal Dimension, }%
Springer-Verlag, New York, 1995.
\end{thebibliography}
\end{document}